\documentclass{llncs}

\usepackage{pstricks,pst-node}
\usepackage{amsmath,amsfonts,amssymb}

\begin{document}

\title{
	Isogenies and the Discrete Logarithm Problem 
	in Jacobians of Genus 3 Hyperelliptic Curves
}
\titlerunning{Isogenies and the DLP on hyperelliptic genus~$3$ curves}

\author{
	Benjamin Smith
}
\institute{ 
	INRIA Saclay--\^Ile-de-France 
	\\
	Laboratoire d'Informatique de l'\'Ecole polytechnique (LIX)
	\\
	91128 Palaiseau Cedex, France
	\email{smith@lix.polytechnique.fr}
	\url{http://www.lix.polytechnique.fr/Labo/Benjamin.Smith/}
}


\newtheorem{hypothesis}{Hypothesis}
\newtheorem{algorithm}[definition]{Algorithm}


\renewcommand{\AA}[2][{}]{\ensuremath{\mathbb{A}^{#2}_{#1}}}
\newcommand{\CC}{{\mathbb{C}}}
\newcommand{\FF}{{\mathbb{F}}}
\newcommand{\FFbar}{{\overline{\mathbb{F}}}}
\newcommand{\QQ}{{\mathbb{Q}}}
\newcommand{\ZZ}{{\mathbb{Z}}}
\newcommand{\NN}{{\mathbb{N}}}
\newcommand{\PP}[2][{}]{\ensuremath{\mathbb{P}^{#2}_{#1}}}
\newcommand{\Jac}[1]{\ensuremath{{J}_{#1}}}
\newcommand{\Pic}[2][{}]{\ensuremath{\mathrm{Pic}^{#1}(#2)}}
\newcommand{\Div}[2][{}]{\ensuremath{\mathrm{Div}^{#1}(#2)}}
\newcommand{\Prin}[1]{\ensuremath{\mathrm{Prin}(#1)}}
\newcommand{\Gal}[1]{{\mathrm{Gal}(#1)}}
\newcommand{\Aut}{{\mathrm{Aut}}}
\newcommand{\Hom}{{\mathrm{Hom}}}
\newcommand{\End}{{\mathrm{End}}}
\newcommand{\Mat}{{\mathrm{Mat}}}
\newcommand{\subgroup}[1]{\ensuremath{\left\langle{#1}\right\rangle}}
\newcommand{\variety}[1]{\ensuremath{V\!\left({#1}\right)}}
\newcommand{\ideal}[1]{\ensuremath{\left(#1\right)}}
\newcommand{\dualof}[1]{\ensuremath{{#1}^\dagger}}
\newcommand{\compose}[2]{\ensuremath{{#1}\circ{#2}}}
\newcommand{\multiplication}[2][{}]{\ensuremath{[#2]_{#1}}}
\newcommand{\genus}[1]{\ensuremath{g_{#1}}}
\newcommand{\differentials}{\ensuremath{\Omega}}
\newcommand{\Spec}{\mathrm{Spec}}
\newcommand{\moduli}[1]{\ensuremath{\mathcal{M}_{#1}}}
\newcommand{\hyperellipticmoduli}[1]{\ensuremath{\mathcal{H}_{#1}}}
\newcommand{\abelianmoduli}[1]{\ensuremath{\mathcal{A}_{#1}}}
\newcommand{\ord}{\mathrm{ord}}
\newcommand{\classof}[1]{\ensuremath{\left[{#1}\right]}}
\newcommand{\integers}[1]{\ensuremath{\family{O}_{#1}}}
\newcommand{\isogenytype}[2]{\ensuremath{(\ZZ/{#1}\ZZ)^{#2}}}
\newcommand{\differentialmatrix}[1]{\ensuremath{M({#1})}}
\newcommand{\Gr}{\mathrm{Gr}}
\newcommand{\softO}{\ensuremath{\widetilde O}}
\newcommand{\ul}[1]{\ensuremath{\underline{#1}}}
\newcommand{\tractables}[1]{\ensuremath{\mathcal{S}}({#1})}
\newcommand{\homog}[1]{\ensuremath{\widetilde{#1}}}
\newcommand{\Galois}{\ensuremath{\mathcal{G}}}


\maketitle

\begin{abstract}
We describe the use of explicit isogenies
to translate instances of the Discrete Logarithm Problem~(DLP)
from Jacobians of hyperelliptic genus~$3$ curves
to Jacobians of non-hyperelliptic genus~$3$ curves,
where they are vulnerable to faster index calculus attacks.
We provide explicit formulae for isogenies
with kernel isomorphic to~$(\ZZ/2\ZZ)^3$
(over an algebraic closure of the base field)
for any hyperelliptic genus~$3$ curve
over a field of characteristic not~$2$ or~$3$.
These isogenies are rational
for a positive fraction of all hyperelliptic genus~$3$ curves
defined over a finite field of characteristic~$p > 3$.
Subject to reasonable assumptions,
our constructions give an explicit and efficient
reduction of instances of the DLP from hyperelliptic 
to non-hyperelliptic Jacobians
for around~$18.57\%$ of all hyperelliptic genus~$3$ curves
over a given finite field.
We conclude with a discussion on extending these ideas
to isogenies with more general kernels.
A condensed version of this work appeared in the proceedings of
the EUROCRYPT 2008 conference.
\end{abstract}

\section{Introduction}
\label{section:introduction}

After the great success of elliptic curves in public-key cryptography,
researchers have naturally been drawn to their higher-dimensional
generalizations: Jacobians of higher-genus curves.
Curves of genus~$1$ (elliptic curves),~$2$, and~$3$ are widely believed
to offer the best balance of security and efficiency.
This article is concerned with the security of curves of genus~$3$.

There are two classes of curves of genus~$3$:
hyperelliptic and non-hyperelliptic.
Each class has a distinct geometry:
the canonical morphism of a hyperelliptic curve
is a double cover of a curve of genus~$0$,
while the canonical morphism of a non-hyperelliptic curve of genus~$3$
is a birational map to a nonsingular plane quartic curve.
A hyperelliptic curve cannot be isomorphic (or birational)
to a non-hyperelliptic curve.
From a cryptological point of view,
the Discrete Logarithm Problem (DLP)
in Jacobians of hyperelliptic curves of genus~$3$ over~$\FF_{q}$
may be solved in~$\softO(q^{4/3})$ group operations,
using the index calculus algorithm of 
Gaudry, Thom\'e, Th\'eriault, and Diem~\cite{GTTD}.
Jacobians of non-hyperelliptic curves of genus~$3$ over~$\FF_{q}$
are amenable to 
Diem's index calculus algorithm~\cite{Diem},
which requires only~$\softO(q)$ 
group operations
to solve the DLP
(for comparison, Pollard/baby-step-giant-step methods 
require~$\softO(q^{3/2})$ group operations
to solve the DLP in Jacobians 
of genus~$3$ curves over~$\FF_{q}$).
The security of non-hyperelliptic genus~$3$ curves
is therefore widely held to be lower than
that of their hyperelliptic cousins.

Our aim 
is to construct explicit homomorphisms
to provide a means of efficiently translating instances of the DLP
from Jacobians of hyperelliptic curves of genus~$3$
to Jacobians of non-hyperelliptic curves,
where faster index calculus is available.
In the context of DLP-based cryptography,
we may assume that our Jacobians are absolutely simple.
In this situation, 
every nontrivial homomorphism of Jacobians of curves of genus~$3$
is an \emph{isogeny}:
that is, a surjective homomorphism with finite kernel.

To be specific,
suppose we are given a hyperelliptic curve~$H$ of genus~$3$
over a finite field~$\FF_{q}$,
together with an instance~$P = [n]Q$ of the DLP in~$\Jac{H}(\FF_{q})$;
our task is to recover~$n$ given~$P$ and~$Q$.
After applying the standard Pohlig--Hellman reduction~\cite{Pohlig--Hellman},
we may assume that~$P$ and~$Q$ have prime order.
We want to solve this DLP instance 
by solving an equivalent DLP instance in a non-hyperelliptic Jacobian.
Suppose we have an isogeny~$\phi: \Jac{H} \to \Jac{C}$,
where~$C$ is a non-hyperelliptic curve of genus~$3$.
Further,
suppose that~$\phi$ is explicit
(that is, we have equations for~$C$
and an efficient map on divisor classes representing~$\phi$)
and defined over~$\FF_{q}$,
so it maps~$\Jac{H}(\FF_{q})$ into~$\Jac{C}(\FF_{q})$.
Provided~$\phi(Q) \not= 0$, 
we can recover~$n$
by solving the DLP instance~$\phi(P) = [n]\phi(Q)$
in~$\Jac{C}(\FF_{q})$ with Diem's algorithm.

The approach outlined above is conceptually straightforward;
the difficulty lies in computing 
explicit isogenies of Jacobians of genus~$3$ curves.
Automorphisms, integer multiplications, and Frobenius maps
aside,
we know of no explicit and general formulae for isogenies from
Jacobians of hyperelliptic curves of genus~$3$
apart from those presented below.

In~\S\ref{section:kernel} through~\S\ref{section:equations-for-the-isogeny}, 
we derive explicit formulae
for isogenies whose kernels are generated by
differences of Weierstrass points,
following the construction of Donagi and Livn\'e~\cite{Donagi--Livne}.
The key step is making Recillas' trigonal construction~\cite{Recillas}
completely explicit.
This gives us a curve~$X$ of genus~$3$
and an explicit isogeny~$\Jac{H} \to \Jac{X}$.
While~$X$ may be hyperelliptic,
na\"ive moduli space dimension arguments suggest 
(and experience confirms) 
that~$X$ will be non-hyperelliptic
with an overwhelming probability,
and thus explicitly isomorphic to a nonsingular plane quartic curve~$C$.
We can therefore compute an explicit isogeny~$\phi: \Jac{H}\to\Jac{C}$;
if~$\phi$ is defined over~$\FF_{q}$,
then we can use it to reduce DLP instances.
We note that the trigonal construction
(and hence our formulae)
does not apply in characteristics~$2$ and~$3$.

We show in \S\ref{section:probabilities}
that, subject to some reasonable assumptions,
given a uniformly randomly chosen hyperelliptic curve~$H$ 
of genus~$3$ over a sufficiently large
finite field~$\FF_{q}$ of characteristic at least~$5$,
our algorithms succeed in constructing an explicit isogeny 
defined over~$\FF_{q}$
from~$\Jac{H}$ to a non-hyperelliptic Jacobian
with probability~$\approx 0.1857$.
In particular,
instances of the DLP 
can be solved in~$\softO(q)$ group operations
for around~$18.57\%$ of all Jacobians of hyperelliptic curves of genus~$3$
over finite fields of characteristic at least~$5$.

We discuss more general isogenies in~\S\ref{section:other-isogenies}.
Given explicit formulae for these isogenies,
we expect that most, if not all, instances of the DLP
in Jacobians of hyperelliptic curves of genus~$3$
over any finite field
could be reduced to instances of the DLP in non-hyperelliptic Jacobians.

Our results have a number of interesting implications
for curve-based cryptography,
at least for curves of genus~$3$.
First,
the difficulty of the DLP in a subgroup~$G$ of~$\Jac{H}$
depends not only on the size of the subgroup~$G$,
but upon the existence of other rational subgroups of~$\Jac{H}$
that can be used to form quotients.
Second,
the security of a given hyperelliptic genus~$3$ curve
depends significantly upon the factorization of its hyperelliptic polynomial.
Neither of these results has any parallel
in genus~$1$ or~$2$.

The constructions of \S\ref{section:kernel}
through \S\ref{section:equations-for-the-isogeny}
and \S\ref{section:other-isogenies}
require some nontrivial algebraic geometry.
We have included enough mathematical detail here
to enable the reader to compute examples,
to justify our claim that the construction is efficient, 
and to support our heuristics.

\subsection*{A Note on the Text}

This article presents an extended version of work that appeared in 
the proceedings of the EUROCRYPT 2008 conference~\cite{Smith-Eurocrypt}.
The chief results are the same;
we have made some (minor) changes to our notation,
expanded the derivation in~\S\ref{section:equations-for-the-isogeny},
given further details and proofs throughout,
and added an appendix
with algorithms to compute sets of tractable subgroups.

\section{Notation and Conventions for Hyperelliptic Curves}
\label{section:notation-and-conventions-for-hyperelliptic-curves}

We will work over~$\FF_{q}$
throughout this article,\footnote{
	Some of the theory 
	carries over to more general base fields:
	in particular,
	the results 
	of~\S\ref{section:computing-trigonal-maps} 
	and~\S\ref{section:equations-for-the-isogeny}
	are valid over fields of characteristic not~$2$ or~$3$.
}
where
$q$ is a power of a prime~$p > 3$.
We let~$\Galois$ denote the Galois group~$\Gal{\FFbar_{q}/\FF_{q}}$,
which is (topologically) generated by the~$q^\mathrm{th}$ power Frobenius map.

Suppose we are given a hyperelliptic curve~$H$ of genus~$3$
over~$\FF_{q}$.
We will use both an affine model
$$
	H: y^2 = F(x) ,
$$
where~$F$ is a squarefree polynomial of degree~$7$ or~$8$,
and a weighted projective plane model
$$
	H: w^2 = \homog{F}(u,v)
$$
for~$H$
(here~$u$,~$v$, and~$w$ 
have weights~$1$,~$1$, and~$4$, respectively).
The coordinates of these models 
are related by 
$x = u/v$ and~$y = w/v^4$.
The polynomial~$\homog{F}$ is squarefree of total degree~$8$,
with
$\homog{F}(u,v) = v^8F(u/v)$ and~$F(x) = \homog{F}(x,1)$.
We emphasize that~$F$ need not be monic.
By a \emph{randomly chosen hyperelliptic curve},
we mean the hyperelliptic curve defined by~$w^2 = \homog{F}(u,v)$,
where~$\homog{F}$ is a uniformly randomly chosen 
squarefree homogenous bivariate polynomial of degree~$8$
over~$\FF_{q}$.

The canonical \emph{hyperelliptic involution}~$\iota$ of~$H$
is defined by 
$(x,y) \mapsto (x,-y)$ in the affine model, 
$(u:v:w) \mapsto (u:v:-w)$ in the projective model,
and induces the negation map~$[-1]$ on~$\Jac{H}$.
The quotient~$\pi: H \to H/\subgroup{\iota}\cong \PP{1}$
sends~$(u:v:w)$ to~$(u:v)$ in the projective model,
and~$(x,y)$ to~$x$ in the affine model 
(where it maps onto the affine patch of~$\PP{1}$ where~$v \not= 0$).

To compute in~$\Jac{H}$,
we fix an isomorphism from~$\Jac{H}$ 
to the group of degree-$0$ divisor classes on~$H$, 
denoted~$\Pic[0]{H}$.
Recall that divisors are formal sums of points in~$H(\FFbar_q)$,
and if
\(
	D = \sum_{P \in H} n_P(P)
\)
is a divisor,
then~$\sum_{P\in H} n_P$ is the \emph{degree} of~$D$.
We say~$D$ is \emph{principal}
if 
\(
	D = \mathrm{div}(f) := \sum_{P\in H} \mathrm{ord}_P(f)(P)
\)
for some function~$f$ on~$H$,
where~$\mathrm{ord}_P(f)$
denotes the number of zeroes (or the negative of the number of poles)
of~$f$ at~$P$.
Since~$H$ is complete,
every principal divisor has degree~$0$.
The group~$\Pic[0]{H}$
is defined to be the group of divisors of degree~$0$
modulo principal divisors;
the equivalence class of a divisor~$D$ is denoted by~$[D]$.
We let~$\Jac{H}[l]$ denote the~$l$-torsion subgroup of~$\Jac{H}$:
that is, the kernel of the multiplication-by-$l$ map.
If~$l$ is prime to~$q$,
then~$\Jac{H}[l](\FFbar_{q})$ is isomorphic to~$(\ZZ/l\ZZ)^6$.

\section{The Kernel of the Isogeny}
\label{section:kernel}

The eight points of~$H(\FFbar_{q})$ where~$w = 0$
are called the \emph{Weierstrass points} of~$H$.
Each Weierstrass point~$W$ 
corresponds to a linear factor
\[
	L_W := v(W)u - u(W)v
\]
of~$\homog{F}$,
which is defined up to scalar multiples.
If~$W_1$ and~$W_2$
are Weierstrass points,
then
\(
	2(W_1) - 2(W_2) 
	= 
	\mathrm{div}(L_{W_1}/L_{W_2})
\),
so~$2[(W_1) - (W_2)] = 0$; hence~$[(W_1) - (W_2)]$
represents an element of~$\Jac{H}[2](\FFbar_{q})$.
In particular,~$[(W_1) - (W_2)] = [(W_2) - (W_1)]$,
so the divisor class~$[(W_1) - (W_2)]$
corresponds to the pair~$\{W_1, W_2\}$
of Weierstrass points,
and hence to the quadratic factor~$L_{W_1}L_{W_2}$ of~$\homog{F}$
(up to scalar multiples).

\begin{proposition}
\label{proposition:partitions-and-tractable-subgroups}
	To every~$\Galois$-stable partition 
	of the eight Weierstrass points of~$H$
	into four disjoint pairs,
	we may associate
	an~$\FF_{q}$-rational subgroup of 
	$\Jac{H}[2](\FFbar_{q})$
	isomorphic to~$(\ZZ/2\ZZ)^3$.
\end{proposition}
\begin{proof}
	Let
	\(
		\{
			\{W_1',W_1''\}, \{W_2',W_2''\},
			\{W_3',W_3''\}, \{W_4',W_4''\}
		\}
	\)
	be a partition of the set of Weierstrass points of~$H$
	into four disjoint pairs.
	Each pair~$\{W_i', W_i''\}$
	corresponds to the~$2$-torsion divisor class
	$[(W_i') - (W_i'')]$ in~$\Jac{H}[2](\FFbar_{q})$.
	We associate the subgroup
	$S := \subgroup{[(W_i') - (W_i'') ] : 1 \le i \le 4}$
	to the partition.
	Observe that
	\[
		\sum_{i = 1}^4 [ (W_i') - (W_i'') ] 
		= 
		\Big[\mathrm{div}\big(w/\prod_{i=1}^4 L_{W_i''}\big)\Big] 
		= 
		0
		;
	\]
	this is the only relation on the classes~$[(W_i') - (W_i'')]$,
	so~$S \cong (\ZZ/2\ZZ)^3~$.
	The action of~$\Galois$
	on~$\Jac{H}[2](\FFbar_{q})$
	corresponds to its action on the Weierstrass points,
	so if the partition is~$\Galois$-stable,
	then the subgroup~$S$ is~$\Galois$-stable.
\qed
\end{proof}

\begin{remark}
	By ``an~$\FF_{q}$-rational subgroup of~$\Jac{H}[2](\FFbar_{q})$
	isomorphic to~$(\ZZ/2\ZZ)^3$'',
	we mean a~$\Galois$-stable subgroup 
	that is isomorphic to~$(\ZZ/2\ZZ)^3$
	over~$\FFbar_{q}$.
	We emphasize that 
	the subgroup 
	need \emph{not} be contained in~$\Jac{H}(\FF_{q})$.
\end{remark}

\begin{remark}
	Requiring the pairs of Weierstrass points
	in Proposition~\ref{proposition:partitions-and-tractable-subgroups}
	to be disjoint ensures that the associated subgroup
	is isotropic with respect to the $2$-Weil pairing.
	We will see in~\S\ref{section:other-isogenies} that
	this is necessary for the quotient by the subgroup
	to be an isogeny of principally polarized abelian varieties,
	and hence for the quotient 
	to be an isogeny of Jacobians.
\end{remark}

\begin{definition}
	We call the subgroups corresponding
	to partitions of the Weierstrass points of~$H$
	as in
	Proposition~\ref{proposition:partitions-and-tractable-subgroups}		\emph{tractable subgroups}.
	We let~$\tractables{H}$ denote
	the set of all~$\FF_{q}$-rational	
	tractable subgroups of~$\Jac{H}[2](\FFbar_{q})$.
\end{definition}

\begin{remark}
	Not every subgroup of~$\Jac{H}[2](\FFbar_{q})$ 
	that is the kernel of an isogeny of Jacobians 
	is a tractable subgroup.
	For example,
	if~$W_1,\ldots,W_8$ are the Weierstrass points of~$H$,
	then the subgroup 
	\[
		\subgroup{
		[ (W_1) - (W_i) + (W_j) - (W_k) ]
		:
			(i,j,k) \in \{ (2,3,4), (2,5,6), (3,5,7) \}
		}
	\]
	is a maximal~$2$-Weil isotropic subgroup of $\Jac{H}(\FFbar_{q})$, 
	and hence is the kernel
	of an isogeny of Jacobians
	(see~\S\ref{section:other-isogenies}).
	However,
	this subgroup contains no nontrivial differences of Weierstrass points,
	and therefore cannot be a tractable subgroup.
\end{remark}	

Computing~$\tractables{H}$ is straightforward
if we identify each tractable subgroup with
its corresponding partition of Weierstrass points.
Recall that each pair of Weierstrass points~$\{W_i',W_i''\}$
corresponds to a quadratic factor of~$\homog{F}$ (up to scalar multiples).
Since the pairs are disjoint,
the corresponding quadratic factors are pairwise coprime,
so we may take them to form 
a factorization of~$\homog{F}$.
We therefore have a correspondence
of tractable subgroups,
partitions of Weierstrass points into pairs,
and sets of quadratic polynomials (up to scalar multiples):
\[
	S 
	\longleftrightarrow 
	\big\{ \{ W_i', W_i'' \} : 1 \le i \le 4 \big\}
	\longleftrightarrow
	\big\{ F_1, F_2, F_3, F_4 \big\}, 
	\ \text{where} \ 	
	\homog{F} = F_1F_2F_3F_4
	.
\]
The action of~$\Galois$
on~$\Jac{H}[2](\FFbar_{q})$
corresponds to its action on the set of Weierstrass points,
so the action of~$\Galois$ 
on a tractable subgroup~$S$
corresponds to the action of~$\Galois$
on the corresponding set~$\{ F_1, F_2, F_3, F_4 \}$
(assuming the $F_i$ have been scaled appropriately).
In particular, 
$S$ is~$\FF_{q}$-rational 
precisely when~$\{ F_1, F_2, F_3, F_4 \}$
is fixed by~$\Galois$.
The factors~$F_i$ are themselves defined over~$\FF_{q}$
precisely when the corresponding points of~$S$ are~$\FF_{q}$-rational.

We can use this information to compute~$\tractables{H}$.
The set of pairs of Weierstrass points
contains a~$\Galois$-orbit
\(
	\big( \{ W_{i_1}',W_{i_1}'' \}, \ldots, \{ W_{i_n}',W_{i_n}'' \} \big) 
\)
if and only if 
(possibly after exchanging some of the~$W_{i_k}'$ with the~$W_{i_k}''$) 
either
both~$( W_{i_1}', \ldots, W_{i_n}' )$
and~$( W_{i_1}'', \ldots, W_{i_n}'' )$
are~$\Galois$-orbits
or
$ ( W_{i_1}', \ldots, W_{i_n}',  W_{i_1}'', \ldots, W_{i_n}'' )~$
is a~$\Galois$-orbit.
Every~$\Galois$-orbit of Weierstrass points
corresponds to an~$\FF_{q}$-irreducible factor of~$F$,
so the size of~$\tractables{H}$
depends only on the factorization of~$F$.
A table 
relating the size of~$\tractables{H}$
to the factorization of~$\homog{F}$
appears in Lemma~\ref{lemma:number-of-tractable-subgroups}
below;
this will be useful
for our analysis in~\S\ref{section:probabilities}.
For completeness,
we have included a na\"ive algorithm
for enumerating~$\tractables{H}$
in Appendix~\ref{appendix:computing-SH}.

\begin{lemma}
\label{lemma:number-of-tractable-subgroups}
	Let~$H: w^2 = \homog{F}(u,v)$ be a hyperelliptic curve of genus~$3$
	over~$\FF_{q}$.
	The cardinality of the set~$\tractables{H}$
	depends only on 
	the degrees of the~$\FF_{q}$-irreducible factors of~$\homog{F}$,
	and is described by the following table:
	\begin{center}
	\begin{tabular}{|c|c|}
		\hline
		Degrees of~$\FF_{q}$-irreducible factors of~$\homog{F}$
		&
		$\# \tractables{H}$
		\\
		\hline
		\hline
		$( 8 ), (6, 2), (6, 1, 1), (4,2,1,1)$ &~$1$ \\ \hline
		$( 4, 2, 2 ), (4,1,1,1,1), (3,3,2), (3,3,1,1)$ &~$3$ \\ \hline
		$( 4, 4 )$ &~$5$ \\ \hline
		$( 2, 2, 2, 1, 1 )$ &~$7$ \\ \hline
		$( 2, 2, 1, 1, 1, 1 )$ &~$9$ \\ \hline
		$( 2, 1, 1, 1, 1, 1, 1 )$ &~$15$ \\ \hline
		$( 2, 2, 2, 2 )$ &~$25$ \\ \hline
		$( 1, 1, 1, 1, 1, 1, 1, 1 )$ &~$105$ \\ \hline
		Other &~$0$ \\ \hline
	\end{tabular}
	\end{center}
\end{lemma}
\begin{proof}
	This is a routine combinatorial exercise
	after noting that every~$\Galois$-orbit
	of pairs of Weierstrass points
	corresponds to either an even-degree factor of~$F$,
	or a pair of factors of~$F$ of the same degree.
\qed
\end{proof}

\section{The Trigonal Construction}
\label{section:trigonal-construction}

We will now briefly outline the theoretical aspects
of constructing isogenies with tractable kernels.
We will make the construction completely explicit in 
\S\ref{section:computing-trigonal-maps}
and~\S\ref{section:equations-for-the-isogeny}.


\begin{definition}
	Suppose~$S = \subgroup{[(W_i') - (W_i'')] : 1 \le i \le 4}$
	is a tractable subgroup.
	We say that a morphism~$g: \PP{1} \to \PP{1}$
	is a \emph{trigonal map for~$S$}
	if 
	$g$ has degree~$3$
	and~$g(\pi(W_i')) = g(\pi(W_i''))$ for~$1 \le i \le 4$.
\end{definition}

Given a trigonal map~$g$ for some tractable subgroup~$S$, 
Recillas' trigonal construction~\cite{Recillas} specifies
a curve~$X$ of genus~$3$ and a map~$f: X \to \PP{1}$ of degree~$4$.\footnote{
	Recillas' original trigonal construction
	is defined where $\pi$ is an \'etale double cover;
	the trigonal construction we apply here 
	is in fact the flat limit of Recillas' construction
	(see~\cite[\S3]{Donagi--Livne} for details).
}
The isomorphism class of~$X$
depends only on~$S$,
and
is independent of the choice of~$g$
(see Recillas~\cite{Recillas}, Donagi~\cite[Th.~2.11]{Donagi}, and Remark~\ref{remark:trigonal-map-descent} below).
Theorem~\ref{theorem:isogeny-from-trigonal-construction},
due to Donagi and Livn\'e,
states that if~$g$ is a trigonal map for~$S$,
then~$S$ is the kernel of an isogeny from~$\Jac{H}$ to~$\Jac{X}$.

\begin{theorem}[Donagi and Livn\'e {\cite[\S5]{Donagi--Livne}}]
\label{theorem:isogeny-from-trigonal-construction}
	Let~$S$ be a tractable subgroup in~$\tractables{H}$,
	and let~$g: \PP{1} \to \PP{1}$
	be a trigonal map for~$S$.
	If~$X$ is the curve formed from~$g$
	by Recillas' trigonal construction,
	then there is an isogeny~$\phi: \Jac{H} \to \Jac{X}$
	(defined over~$\FFbar_{q}$) with kernel~$S$.
\end{theorem}

We will give only a brief description of the geometry of~$X$ here,
concentrating instead on its explicit construction;
we refer the reader to
Recillas~\cite{Recillas}, 
Vakil~\cite{Vakil},
Donagi~\cite[\S2]{Donagi},
and Birkenhake and Lange~\cite[\S12.7]{Birkenhake--Lange}
for proofs and further detail.
The isogeny of Theorem~\ref{theorem:isogeny-from-trigonal-construction}
is analogous to 
the well-known Richelot isogeny
in genus~$2$
(see 
Bost and Mestre~\cite{Bost-Mestre},
and Donagi and Livn\'e~\cite[\S4]{Donagi--Livne}
for details),
and to the explicit isogeny 
described by Lehavi and Ritzenthaler in~\cite{Lehavi--Ritzenthaler}
for Jacobians of non-hyperelliptic genus $3$ curves.

In abstract terms,
if~$U$ is the subset of the codomain of~$g$
above which~$g\circ\pi$ is unramified,
then~$X$ is by definition the closure of
the curve over~$U$
representing the pushforward to~$U$
of the sheaf 
of sections 
of~$\pi: (g\circ\pi)^{-1}(U) \to g^{-1}(U)$
(in the \'etale topology).
This means in particular
that the~$\FFbar_{q}$-points of~$X$ 
over an~$\FFbar_{q}$-point~$P$ of~$U$
represent partitions of 
the six~$\FFbar_{q}$-points of~$(g\circ\pi)^{-1}(P)$
into two sets of three exchanged by the hyperelliptic involution.
The fibre product of~$H$ and~$X$ over~$\PP{1}$
with respect to~$g\circ\pi$ and~$f$
is the union of two isomorphic curves,~$R$ and~$R'$,
which are exchanged by the involution on~$H\times_{\PP{1}}X$
induced by the hyperelliptic involution.
The natural projections induce
coverings~$\pi_H: R \to H$ and~$\pi_X: R \to X$
of degrees~$2$ and~$3$, respectively,
so~$R$ is a~$(3,2)$-correspondence between~$H$ and~$X$.

\begin{figure}
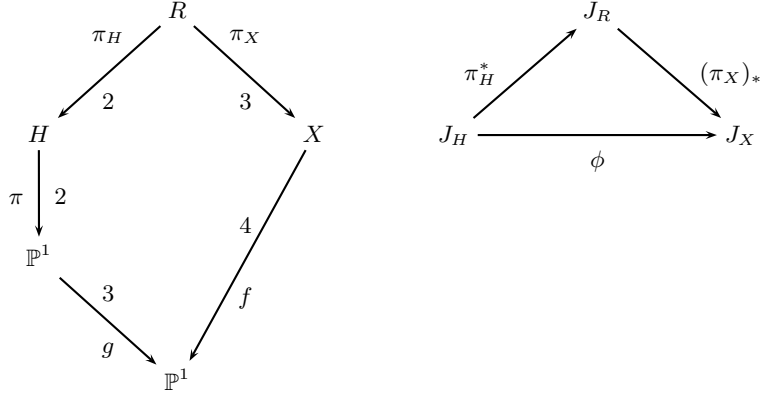

$$
\begin{psmatrix}[rowsep=36pt]
&R&&&\Jac{R}\\
H&&X&\Jac{H}&&\Jac{X}\\
\PP{1}\\
&\PP{1}
\ncline[nodesep=3pt]{->}{1,2}{2,1}^{\pi_H}_2
\ncline[nodesep=3pt]{->}{1,2}{2,3}^{\pi_X}_3
\ncline[nodesep=3pt]{->}{2,1}{3,1}<{\pi}>2
\ncline[nodesep=3pt]{->}{3,1}{4,2}^{3}_g
\ncline[nodesep=3pt]{->}{2,3}{4,2}^{4}_f
\ncline[nodesep=3pt]{->}{2,4}{2,6}_{\phi}
\ncline[nodesep=3pt]{->}{2,4}{1,5}<{\pi_H^*}
\ncline[nodesep=3pt]{->}{1,5}{2,6}>{(\pi_X)_*}
\end{psmatrix}
$$
\caption{The curves, Jacobians, and morphisms of \S\ref{section:trigonal-construction}}
\label{table:curve-diagram}
\end{figure}

The maps~$\pi_H$ and~$\pi_X$
induce homomorphisms~$(\pi_H)^*: \Jac{H} \to \Jac{R}$
(the \emph{pullback})
and~$(\pi_X)_*: \Jac{R} \to \Jac{X}$ (the \emph{pushforward}).
In terms of divisor classes,
the pullback is defined by
\[
	(\pi_H)^*\Big(\Big[\sum_{P \in H}n_P(P)\Big]\Big) 
	= 
	\Big[ \sum_{P \in H}n_P\!\!\!\!{\sum_{Q \in \pi_H^{-1}(P)}}\!\!\!\!(Q) \Big],
\]
with appropriate multiplicities where~$\pi_H$ ramifies;
the pushforward is defined by
\[
	(\pi_X)_*\Big(\Big[ \sum_{Q\in R}m_Q(Q) \Big]\Big) 
	= 
	\Big[ \sum_{Q\in R}m_Q(\pi_X(Q)) \Big]
	.
\]
Composing~$(\pi_X)_*$ with~$(\pi_H)^*$,
we obtain
an isogeny~$\phi: \Jac{H} \to \Jac{X}$
with kernel~$S$.
If we replace~$R$ with~$R'$ in the above,
we obtain an isogeny isomorphic to~$-\phi$.
Thus, up to isomorphism, the construction of the isogeny
depends only on the subgroup~$S$.
The curves and Jacobians described above form 
the commutative diagrams shown in~Figure~\ref{table:curve-diagram}.

The hyperelliptic Jacobians form a codimension-$1$ subspace $\hyperellipticmoduli{g}$
of the moduli space of 
$3$-dimensional principally polarized abelian varieties
--- which, by the theorem of Oort and Ueno~\cite{Oort--Ueno},
is also the moduli space $\moduli{g}$ of Jacobians of genus~$3$ curves.
The Weil hypotheses imply that 
$\#\hyperellipticmoduli{g}(\FF_{q})/\#\moduli{g}(\FF_{q}) \sim 1/q$
for sufficiently large $q$
(cf.~\cite[Theorem 1]{Lang-Weil}).
In particular, 
for cryptographically relevant sizes of~$q$,
the probability that a uniformly randomly chosen curve $X$ of genus $3$
over $\FF_{q}$
should be hyperelliptic is negligible.
We will suppose that
the same is true for the curve $X$
constructed in Theorem~\ref{theorem:isogeny-from-trigonal-construction}
for a uniformly randomly chosen $H$ and $S$ in $\tractables{H}$.
This is consistent with our experimental observations,
so we postulate Hypothesis~\ref{hypothesis:hyperellipticity}.

\begin{hypothesis}
\label{hypothesis:hyperellipticity}
	The probability that 
	the curve~$X$
	constructed by the trigonal construction
	for a randomly chosen~$H/\FF_{q}$ and~$S$ in~$\tractables{H}$
	is hyperelliptic is negligible
	for sufficiently large $q$.
\end{hypothesis}

\section{Computing Trigonal Maps}
\label{section:computing-trigonal-maps}

Suppose we are given a tractable subgroup
$S$ of~$\Jac{H}[2](\FFbar_{q})$,
corresponding to a partition~$\{ \{W_i', W_i''\}: 1 \le i \le 4 \}$
of the Weierstrass points of~$H$ into pairs.
The first step in the explicit trigonal construction
is to compute a trigonal map~$g$ for~$S$.
We will
compute polynomials~$N = x^3 + n_1x + n_0$ 
and~$D = x^2 + d_1x + d_0$
such that the rational map 
\begin{equation}
\label{eq:rational-map-form}
	g: x \longmapsto t 
	= 
	\frac{N(x)}{D(x)} 
	= 
	\frac{x^3 + n_1x + n_0}{x^2 + d_1x + d_0}
\end{equation}
defines a trigonal map for~$S$.
The derivation is an exercise in classical geometry;
we include it here to demonstrate its efficiency
and to justify Hypothesis~\ref{hypothesis:trigonal-map-rationality},
which will be important in determining the expectation of success
of our reduction in~\S\ref{section:probabilities}.
The reader prepared to admit
the existence of efficiently computable trigonal maps 
in the form of~\eqref{eq:rational-map-form} 
may skip the remainder of this section on first reading.


By definition,
$g: \PP{1} \to \PP{1}$ is
a degree-$3$ map
with
$ g(\pi(W_i')) = g(\pi(W_i''))$ for~$1 \le i \le 4$.
We will express~$g$ as a composition 
$ g = p\circ e$,
where~$e: \PP{1} \to \PP{3}$ is the rational normal embedding
defined by
\[
	e: (u:v) \longmapsto (u_0:u_1:u_2:u_3) = (u^3 : u^2v : uv^2 : v^3) ,
\]
and~$p: \PP{3} \to \PP{1}$ is the projection defined as follows.
For each~$1 \le i \le 4$,
we let~$L_i$ denote the line in~$\PP{3}$
passing through~$e(\pi(W_i'))$ and~$e(\pi(W_i''))$.
There exists at least one line~$L$
intersecting all four of the~$L_i$
(in fact there are two, though they may coincide; we will compute them below).
We take~$p$
to be the projection away from~$L$;
then 
$p(e(\pi(W_i'))) = p(e(\pi(W_i'')))$
for~$1 \le i \le 4$,
so
$g = p\circ e$ is a trigonal map for~$S$.
Given linear equations for~$L$ in the coordinates~$u_i$,
we can use Gaussian elimination to
compute elements~$n_1$,~$n_0$,~$d_1$, and~$d_0$ of~$\FF_{q}$
such that
\[
	L = \variety{ u_0 + n_1u_2 + n_0u_3, u_1 + d_1u_2 + d_0u_3 } .
\]
The projection~$p: \PP{3} \to \PP{1}$ away from~$L$
is then defined by
\[
	p: (u_0:u_1:u_2:u_3)
	\longmapsto
	( u_0 + n_1u_2 + n_0u_3 : u_1 + d_1u_2 + d_0u_3 )
	,
\]
so our trigonal map~$g = p\circ e$
is defined by
\[
	g: (u:v)
	\longmapsto
	( u^3 + n_1uv^2 + n_0v^3 : u^2v + d_1uv^2 + d_0v^3 )
	.
\]
Therefore,
if we set~$N(x) := x^3 + n_1x + n_0$ and~$D(x) := x^2 + d_1x + d_0$,
then~$g$ will be defined by the rational map
$x \longmapsto t = N(x)/D(x)$.

To compute equations for~$L$,
we will use the classical theory of \emph{Grassmannian varieties}.
The elementary Lemmas~\ref{lemma:Grassmannian-correspondence}
and~\ref{lemma:hyperplane}
will be stated without proof;
we refer the reader to
Griffiths and Harris~\cite[\S1.5]{Griffiths--Harris}
and Harris~\cite[Lecture 6]{Harris}
for details.

The set of lines in~$\PP{3}$
has the structure of an algebraic variety~$\Gr(1,3)$,
called the Grassmannian.
There is a convenient model for~$\Gr(1,3)$ 
as a quadric hypersurface in~$\PP{5}$:
if~$v_0,\ldots,v_5$ are coordinates on~$\PP{5}$,
then we may take
\[
	\Gr(1,3) := \variety{ v_0v_3 + v_1v_4 + v_2v_5 } 
	\subset \PP{5}.
\]

\begin{lemma}
\label{lemma:Grassmannian-correspondence}
	There is a bijection between points of~$\Gr(1,3)(\FFbar_{q})$
	and lines in~$\PP{3}$,
	defined as follows.
	\begin{enumerate}
		\item
			The point of~$\Gr(1,3)(\FFbar_{q})$
			corresponding to the line
			through 
			$(p_0:p_1:p_2:p_3)$ and~$(q_0:q_1:q_2:q_3)$ 
			in~$\PP{3}$
			has coordinates
			\[
				\left(
				\left|\begin{array}{@{\:}c@{\;\;}c@{\:}}
					p_0 & p_1 \\
					q_0 & q_1
				\end{array}\right|:
				\left|\begin{array}{@{\:}c@{\;\;}c@{\:}}
					p_0 & p_2 \\
					q_0 & q_2
				\end{array}\right|:
				\left|\begin{array}{@{\:}c@{\;\;}c@{\:}}
					p_0 & p_3 \\
					q_0 & q_3
				\end{array}\right|:
				\left|\begin{array}{@{\:}c@{\;\;}c@{\:}}
					p_2 & p_3 \\
					q_2 & q_3
				\end{array}\right|:
				\left|\begin{array}{@{\:}c@{\;\;}c@{\:}}
					p_3 & p_1 \\
					q_3 & q_1
				\end{array}\right|:
				\left|\begin{array}{@{\:}c@{\;\;}c@{\:}}
					p_1 & p_2 \\
					q_1 & q_2
				\end{array}\right|
				\right)
				.
			\] 
		\item	
			The line in~$\PP{3}$ corresponding 
			to a point
			$(\gamma_0:\cdots:\gamma_5)$
			of~$\Gr(1,3)(\FFbar_{q})$
			is defined~by
			\[
				\variety{
    					\begin{array}{r@{\;}r@{\;}r@{\;}r@{\;}r@{\;}r@{\;}r} 
					0 u_0 & - & \gamma_3 u_1 & - & \gamma_4 u_2 & - & \gamma_5 u_3 , \\
					\gamma_3 u_0 & + & 0 u_1 & - & \gamma_2 u_2 & + & \gamma_1 u_3 , \\ 
					\gamma_4 u_0 & + & \gamma_2 u_1 & + & 0 u_2 & - & \gamma_0 u_3 , \\
					\gamma_5 u_0 & - & \gamma_1 u_1 & + & \gamma_0 u_2 & + & 0 u_3 \\ 
					\end{array} 
				}
			\]
			(two of the equations will be 
			redundant linear combinations of the others).
	\end{enumerate}
\end{lemma}

\begin{lemma}
\label{lemma:hyperplane}
	Let~$L$ be the line in~$\PP{3}$
	corresponding to a point 
	$(\gamma_0:\cdots:\gamma_5)$ of~$\Gr(1,3)(\FFbar_{q})$.
	The points in~$\Gr(1,3)(\FFbar_{q})$
	corresponding to lines in~$\PP{3}$ 
	that intersect nontrivially  with~$L$
	are precisely the points lying in the hyperplane
	defined by \( \sum_{i=0}^5 \gamma_iv_{i+3} = 0 \)
	(where the subscripts are taken modulo~$6$).
\end{lemma}

Suppose~$S$
is represented by a set~$\{F_i=a_iu^2 + b_iuv + c_iv^2 : 1\:\le\:i\:\le~4\}$
of quadratic factors of~$\homog{F}$ (as in~\S\ref{section:kernel}),
with each factor $F_i$ corresponding to a pair~$\{ W_i',W_i'' \}$
of Weierstrass points.
Applying Lemma~\ref{lemma:Grassmannian-correspondence},
we see
that
the line~$L_i$ through~$e(\pi(W_i'))$ and~$e(\pi(W_i''))$
corresponds to the point 
\[
		(c_i^2:-c_ib_i:b_i^2-a_ic_i:a_i^2:a_ib_i:a_ic_i) 
\] 
on~$\Gr(1,3)$.
If~$(\gamma_0:\cdots:\gamma_5)$ in~$\Gr(1,3)(\FFbar_{q})$
corresponds to a candidate for~$L$,
then by 
Lemma~\ref{lemma:hyperplane}
we have~$M(\gamma_0,\ldots,\gamma_5)^T = 0$,
where 
\begin{equation}
\label{eq:M-def}
	M = \left(\!\begin{array}{c@{\ \ }c@{\ \ }c@{\ \ }c@{\ \ }c@{\ \ }c}
		a_1^2 & a_1b_1 & a_1c_1 & c_1^2 & - c_1b_1 & (b_1^2-a_1c_1) \\
		a_2^2 & a_2b_2 & a_2c_2 & c_2^2 & - c_2b_2 & (b_2^2-a_2c_2) \\
		a_3^2 & a_3b_3 & a_3c_3 & c_3^2 & - c_3b_3 & (b_3^2-a_3c_3) \\
		a_4^2 & a_4b_4 & a_4c_4 & c_4^2 & - c_4b_4 & (b_4^2-a_4c_4) \\
	\end{array}
	\!\right) .
\end{equation}
The kernel of~$M$ is two-dimensional, 
corresponding to a line~$\Lambda$ in~$\PP{5}$.
The kernel 
is independent of the ordering of the~$F_i$,
and does not change if we replace the~$F_i$
by scalar multiples;
hence,~$\Lambda$ depends only on the subgroup~$S$.
Let
$\{\ul{\alpha},\ul{\beta}\}$
be a basis for~$\ker M$,
writing~$\ul{\alpha} = (\alpha_0,\ldots,\alpha_5)$
and~$\ul{\beta} = (\beta_0,\ldots,\beta_5)$.
If~$S$ is~$\FF_{q}$-rational,
then so is~$\ker M$,
so we may take the~$\alpha_i$ and~$\beta_i$
to be in~$\FF_{q}$
(see Cartier~\cite[\S I]{Cartier}).
We want to find a point 
$P_L = (\alpha_0 + \lambda \beta_0:\cdots:\alpha_5+\lambda\beta_5)$
where~$\Lambda$
intersects with~$\Gr(1,3)$.
The points~$(u_0:\ldots:u_3)$ 
on the line~$L$ in~$\PP{3}$ corresponding to~$P_L$
satisfy
\(
	(M_{\ul{\alpha}} + \lambda M_{\ul{\beta}})
	(u_0,\ldots,u_3)^T = 0 
\),
where
\[
	M_{\ul{\alpha}}
	:=
	\left(\!\begin{array}{r@{\ \ }r@{\ \ }r@{\ \ }r}
	0 & -\alpha_3 & -\alpha_4 & -\alpha_5 \\
	\alpha_3 & 0 & - \alpha_2 & \alpha_1 \\ 
	\alpha_4 & \alpha_2 & 0 & - \alpha_0 \\
	\alpha_5 & -\alpha_1 & \alpha_0 & 0 \\ 
	\end{array}\!\right)
	\ \  \mbox{and}\quad 
	M_{\ul{\beta}}
	:=
	\left(\!\begin{array}{r@{\ \ }r@{\ \ }r@{\ \ }r}
	0 & -\beta_3 & -\beta_4 & -\beta_5 \\
	\beta_3 & 0 & - \beta_2 & \beta_1 \\ 
	\beta_4 & \beta_2 & 0 & - \beta_0 \\
	\beta_5 & -\beta_1 & \beta_0 & 0 \\ 
	\end{array}\!\right)
	.
\]
By part~(2) of Lemma~\ref{lemma:Grassmannian-correspondence},
the rank of~$M_{\ul{\alpha}} + \lambda M_{\ul{\beta}}$
is~$2$.
Using the expression
\begin{equation}
\label{eq:lambda-eqn}
	\det(M_{\ul{\alpha}} + \lambda M_{\ul{\beta}})
	= 
	\Big(
		\frac{1}{2}\big(\sum_{i=0}^6\beta_i\beta_{i+3} \big)\lambda^2
		+ \big(\sum_{i=0}^6 \alpha_i\beta_{i+3} \big)\lambda
		+ \frac{1}{2}\sum_{i=0}^6 \alpha_i\alpha_{i+3} 
	\Big)^2
\end{equation}
(where the subscripts are taken modulo~$6$),
we see that~$M_{\ul{\alpha}} + \lambda M_{\ul{\beta}}$
has rank~$2$ precisely when 
$\det(M_{\ul{\alpha}} + \lambda M_{\ul{\beta}}) = 0$:
we can therefore solve
$\det(M_{\ul{\alpha}} + \lambda M_{\ul{\beta}}) = 0$
to determine a value for~$\lambda$.
Finally, 
we use Gaussian elimination 
to compute~$n_1$,~$n_0$,~$d_1$, and~$d_0$ in~$\FF_{q}(\lambda)$
such that~$(1,0,n_1,n_0)$ and~$(0,1,d_1,d_0)$
generate the rowspace of 
$M_{\underline{\alpha}} + \lambda M_{\underline{\beta}}$.
We then take 
\(
	L = \variety{u_0 + n_1u_2 + n_0u_3, u_1 + d_1u_2 + d_0u_3 } 
\),
and compute~$p$, $e$,
and the trigonal map~$g = p\circ e$ as above.

Since~$L$ is defined over~$\FF_{q}(\lambda)$,
so is the projection~$p$ and the trigonal map~$g$.
But~$\lambda$
satisfies a quadratic equation with coefficients in~$\FF_{q}$,
so
$\FF_{q}(\lambda)$ 
is at most 
a quadratic extension of~$\FF_{q}$.
Computing the discriminant of~$\det(M_{\underline{\alpha}} + \lambda M_{\underline{\beta}})$,
we obtain a criterion for existence of
trigonal maps over~$\FF_{q}$
for a given tractable subgroup.

\begin{proposition}
\label{proposition:trigonal-map-rationality-criterion}
	Suppose~$S$ is a tractable subgroup,
	and let 
	\( \{ 
		\underline{\alpha} = (\alpha_i), 
		\underline{\beta} = (\beta_i)
	\} \)
	be any~$\FF_{q}$-rational basis of the nullspace of the matrix~$M$
	defined in \eqref{eq:M-def}.
	There exists an~$\FF_{q}$-rational
	trigonal map for~$S$
	if and only if
	\begin{equation} 
	\label{eq:discriminant}
		\Big(\sum_{i=0}^6\alpha_i\beta_{i+3}\Big)^2 
		- 
		\Big(\sum_{i=0}^6\alpha_i\alpha_{i+3}\Big)
		\Big(\sum_{i=0}^6\beta_i\beta_{i+3}\Big)
	\end{equation}
	is a square in~$\FF_{q}$,
	where the subscripts are taken modulo~$6$.
\end{proposition}
\begin{proof}
	From the derivation above,
	we see that there exists an $\FF_{q}$-rational trigonal map for $S$
	if and only if we can find a $\lambda$ in $\FF_{q}$
	such that 
	$\det(M_{\underline{\alpha}} + \lambda M_{\underline{\beta}}) = 0$.
	By Equation~\eqref{eq:lambda-eqn},
	we can find such a $\lambda$
	if and only if
	the quadratic polynomial
	\[
		\frac{1}{2}\big(\sum_{i=0}^6\beta_i\beta_{i+3} \big)T^2
		+ \big(\sum_{i=0}^6 \alpha_i\beta_{i+3} \big)T
		+ \frac{1}{2}\sum_{i=0}^6 \alpha_i\alpha_{i+3} 
	\]
	has two roots in $\FF_{q}$.
	This occurs precisely when the discriminant of this polynomial
	--- the expression in~\eqref{eq:discriminant} above ---
	is a square in $\FF_{q}$.
\qed
\end{proof}

Proposition~\ref{proposition:trigonal-map-rationality-criterion}
shows that the rationality of a trigonal map
for a tractable subgroup~$S$
depends only upon whether an element of~$\FF_{q}$
depending only on~$S$
is a square.
It seems reasonable to assume that these field elements
are uniformly distributed for uniformly random choices of~$H$ and~$S$,
and indeed this is consistent with our experimental observations.
Since 
a uniformly randomly chosen element of~$\FF_{q}$
is a square
with probability~$\sim 1/2$,
we propose Hypothesis~\ref{hypothesis:trigonal-map-rationality}.

\begin{hypothesis}
\label{hypothesis:trigonal-map-rationality}
	The probability that 
	there exists an~$\FF_{q}$-rational trigonal map
	for a subgroup $S$ uniformly randomly chosen from $\tractables{H}$,
	where $H$ is a randomly chosen hyperelliptic curve over~$\FF_{q}$,
	is~$1/2$.
\end{hypothesis}

\section{Equations for the Isogeny}
\label{section:equations-for-the-isogeny}

Suppose we have a hyperelliptic curve~$H$
of genus~$3$, a tractable subgroup~$S$ in~$\tractables{H}$,
and a trigonal map~$g$ for~$S$.
We will now perform an explicit trigonal construction on~$g$
to compute a curve~$X$ and an isogeny~$\phi: \Jac{H} \to \Jac{X}$
with kernel~$S$.

We assume that~$g$ has been derived as
in \S\ref{section:computing-trigonal-maps},
and in particular that~$g:\PP{1}\to\PP{1}$ 
is defined by a rational map in the form
\[
	g: 
	x \longmapsto t 
	= \frac{N(x)}{D(x)} 
	= \frac{x^3 + n_1x + n_0}{x^2 + d_1x + d_0}
	.
\]
Observe that~$g$ maps the point at infinity to the point at infinity
(that is,~$(1:0)$).
For notational convenience, we define
\[
	G(t,x) = x^3 + g_2(t)x^2 + g_1(t)x + g_0(t) := N(x) - tD(x) ;
\]
unless otherwise noted, we will view~$G(t,x)$ as an element of~$\FF_{q}[t][x]$.
We have
\[
	g_2(t) = -t,
	\quad
	g_1(t) = n_1 - d_1t, 
	\quad \text{and}\quad  
	g_0(t) = n_0 - d_0t .
\]
We also define~$f_0$,~$f_1$, and~$f_2$ to be the elements of~$\FF_{q}[t]$
such that
\[
	f_0(t) + f_1(t)x + f_2(t)x^2 \equiv F(x) \pmod{ G(t,x) }.
\]

Let~$U$ be the subset of~$\AA{1} = \PP{1}\setminus\{(1:0)\}$
above which~$g\circ\pi$ is unramified.
With the notation above,
\[
	U = \Spec(k[t])\setminus\variety{
		(f_1^2 - 4f_2f_0)(4g_2^3g_0 - g_2^2g_1^2 - 18g_2g_1g_0 + 4g_1^3 + 27g_0^2)
	}
	.
\]
We will derive equations for an affine model~$X|_U$ of~$f^{-1}(U)$
--- that is, the open subset of~$X$ over~$U$.
We will not prove here
that the normalization of~$X|_U$
is isomorphic to the curve~$X$
specified by Recillas,
but we will exhibit a bijection on geometric points.
If~$X$ is not hyperelliptic,
then taking the canonical map of~$X|_U$ into~$\PP{2}$
will give us a nonsingular plane quartic curve~$C$
isomorphic to~$X$.

By definition,
every point~$P$ in~$X|_U(\FFbar_{q})$
corresponds to a pair of unordered triples of points in~$H(\FFbar_{q})$,
exchanged by the hyperelliptic involution,
with each triple
supported on the fibre of~$g\circ\pi$ over~$f(P)$.
To be more explicit,
suppose~$Q$ is a generic point of~$U$.
Since~$g\circ\pi$
is unramified above~$Q$,
we may choose three preimages 
$P_1$,~$P_2$, and~$P_3$ of~$Q$
such that
\[
	(g\circ\pi)^{-1}(Q)
	= 
	\{ P_1, P_2, P_3, \iota(P_1), \iota(P_2), \iota(P_3) \} 
	.
\]
Viewing unordered triples of points
as effective divisors of degree~$3$
(that is, as formal sums of three points),
we have
\begin{equation}
\label{eq:preimage-of-Q-on-X}
	f^{-1}(Q)
	=
	\left\{
	\begin{array}{c}
		Q_1 \leftrightarrow 
		\big\{ 
			P_1 + P_2 + P_3 ,\
			\iota(P_1) + \iota(P_2) + \iota(P_3) 
		\big\}\!, \\
		Q_2 \leftrightarrow
		\big\{ 
			P_1 + \iota(P_2) + \iota(P_3) ,\ 
			\iota(P_1) + P_2 + P_3
		\big\}\!, \\
		Q_3 \leftrightarrow
		\big\{ 
			\iota(P_1) + P_2 + \iota(P_3) ,\ 
			P_1 + \iota(P_2) + P_3 
		\big\}\!, \\
		Q_4 \leftrightarrow
		\big\{ 
			\iota(P_1) + \iota(P_2) + P_3 ,\ 
			P_1 + P_2 + \iota(P_3) 
		\big\}\!\,
	\end{array}
	\right\}
	.
\end{equation}
Note that~$P_i$ and~$\iota(P_i)$
never appear in the same divisor
for any~$1 \le i \le 3$.
There is a one-to-one correspondence
between effective divisors of degree~$3$ on~$H$
satisfying this condition,
and ideals
$(a(x), y - b(x))$
where~$a$ is a monic cubic polynomial
and
$b$ is a quadratic polynomial satisfying~$b^2 \equiv F \pmod{a}$
(this is the well-known \emph{Mumford representation}~\cite[\S IIIa]{Mumford}). 
For example,
$P_1+P_2+P_3$
corresponds to
the ideal~$(a(x),y-b(x))$
where~$a(x) = \prod_i(x - x(P_i))$
and~$b$ satisfies~$y(P_i)~=~b(x(P_i))$
for~$1 \le i \le 3$
(with appropriate multiplicities);
we may compute~$b$ using
the Lagrange interpolation formula.
A divisor is defined over~$\FF_q$
if and only if~$a$ and~$b$ are defined over~$\FF_q$.
The ideal $(a(x),y-b(x))$ corresponds to~$P_1 + P_2 + P_3$
if and only if
$(a(x),y+b(x))$
corresponds to~$\iota(P_1) + \iota(P_2) + \iota(P_3)$;
so each point of~$X$ over~$U$
corresponds to a pair~$\{ (a(x),y\pm b(x)) \}$
of ideals.
We will construct a curve parametrizing these pairs of ideals,
and take this as a model for~$X|_U$.

Suppose~$\{(a(x),y\pm b(x))\}$
is a pair of ideals corresponding to one of the preimages of~$Q$ on~$X|_U$.
The product of the two ideals is equal to the principal ideal~$(a(x))$;
but products of ideals correspond to sums of divisors,
so~$(a(x))$ must cut out the divisor
$P_1 + P_2 + P_3 + \iota(P_1) + \iota(P_2) + \iota(P_3)$ on~$H$.
This divisor is just~$(g\circ\pi)^*(Q)$,
which we know is cut out by~$(G(t(Q),x))$;
so we conclude that
$a(x) = G(t(Q),x)$
for every pair of ideals~$\{(a(x),y\pm b(x))\}$
corresponding to a point in~$f^{-1}(Q)$.
In particular,
the generic point of~$X|_U$
corresponds to a pair of ideals of the form 
$ \{ ( G(t,x), y \pm (b_0 + b_1x + b_2x^2 ) ) \}$,
where~$b_0$,~$b_1$, and~$b_2$
are algebraic functions of~$t$
such that
\begin{equation}
\label{eq:mumford-relation}
	(b_0 + b_1x + b_2x^2)^2 \equiv F(x)
	\pmod{ G(t,x) }
	.
\end{equation}
Viewing~$b_0$,~$b_1$, and~$b_2$
as coordinates on~$\AA{3}$ (over~$\FF_{q}$),
we
expand both sides of~\eqref{eq:mumford-relation}
modulo~$G(t,x)$
and equate coefficients
to obtain a variety~${\widetilde{X}}$ 
in~$U\times\AA{3}$
parametrizing ideals:
\[
	{\widetilde{X}} = \variety{
		\widetilde{c}_0(t,b_0,b_1,b_2),
		\widetilde{c}_1(t,b_0,b_1,b_2),
		\widetilde{c}_2(t,b_0,b_1,b_2)
	}
	,
\]
where
\[
	\begin{array}{r@{\;=\;}l}
	\widetilde{c}_0(t,b_0,b_1,b_2) 
	& 
	g_2(t)g_0(t)b_2^2 - 2g_0(t)b_2b_1 + b_0^2 - f_0(t) 
	,
	\\
 	\widetilde{c}_1(t,b_0,b_1,b_2) 
	& 
	(g_2(t)g_1(t) - g_0(t))b_2^2 - 2g_1(t)b_2b_1 + 2b_1b_0 - f_1(t)
	,\quad \text{and}
	\\
 	\widetilde{c}_2(t,b_0,b_1,b_2) 
	& 
	(g_2(t)^2 - g_1(t))b_2^2 - 2g_2(t)b_2b_1 + 2b_2b_0 + b_1^2 - f_2(t)
	.
	\end{array}
\]

The ideals in each pair~$\{ (G(t,x),y\pm (b_2x^2 + b_1x + b_0)) \}$
are exchanged by the involution
$ \iota_*: {\widetilde{X}} \longrightarrow {\widetilde{X}}~$
defined by
\[
	\iota_*: (t,b_0,b_1,b_2) \longmapsto (t,-b_0,-b_1,-b_2) ;
\]
the curve~$X|_U$ is therefore the quotient of~${\widetilde{X}}$
by~$\subgroup{\iota_*}$.
To make this quotient explicit,
let~$m: U\times\AA{3}\longrightarrow U\times\AA{6}$
be the map defined by
\[
	m : 
	(t,b_0,b_1,b_2) 
	\longmapsto 
	(t,b_{00},b_{01},b_{02},b_{11},b_{12},b_{22})
	=
	(t,b_0^2,b_0b_1,b_0b_2,b_1^2,b_1b_2,b_2^2)
	;
\]
observe that
\[
	m(U\times\AA{3}) = \variety{
		\begin{array}{c}
			b_{01}^2 - b_{00}b_{11} ,\ 
			b_{01}b_{02} - b_{00}b_{12} ,\ 
			b_{02}^2 - b_{00}b_{22} ,\ 
			\\
			b_{02}b_{11} - b_{01}b_{12} ,\ 
			b_{02}b_{12} - b_{01}b_{22} ,\ 
			b_{12}^2 - b_{11}b_{22} 
		\end{array}
	}
	\subset
	U\times \AA{6}
	.
\]
We have~$X|_U = m({\widetilde{X}})$,
so
\[
	X|_U =
	\variety{
	\begin{array}{c}
		c_0(t,b_{00},\ldots,b_{22}),
		c_1(t,b_{00},\ldots,b_{22}),
		c_2(t,b_{00},\ldots,b_{22}),
		\\
		b_{01}^2 - b_{00}b_{11} ,\ 
		b_{01}b_{02} - b_{00}b_{12} ,\ 
		b_{02}^2 - b_{00}b_{22} ,\ 
		\\
		b_{02}b_{11} - b_{01}b_{12} ,\ 
		b_{02}b_{12} - b_{01}b_{22} ,\ 
		b_{12}^2 - b_{11}b_{22} 
		\end{array}
	}
	\subset U\times\AA{6},
\]
where~$c_0$,~$c_1$, and~$c_2$ are the polynomials defined by
\[
	\begin{array}{l}
		c_0(t,b_{00},b_{01},b_{02},b_{11},b_{12},b_{22})
		:=
		g_2g_0b_{22} - 2g_0b_{12} + b_{00} - f_0 ,
		\\
		c_1(t,b_{00},b_{01},b_{02},b_{11},b_{12},b_{22})
		:=
		(g_2g_1 - g_0)b_{22} - 2g_1b_{12} + 2b_{01} - f_1 ,
		\ \text{ and }
		\\
		c_2(t,b_{00},b_{01},b_{02},b_{11},b_{12},b_{22})
		:=
		(g_2^2 - g_1)b_{22} - 2g_2b_{12} + 2b_{02} + b_{11} - f_2 
		.
	\end{array}
\]
Observe that~$X|_U$ is defined over the field of definition of~$g$.

It remains to derive a correspondence~$R$
between~$H$ and~$X|_U$
inducing the isogeny~$\phi$. 
We know that~$R$ is a component of the fibre product
$H\!\times_{\PP{1}}\! X$
(with respect to~$\compose{g}{\pi}$
and~$f$).
We may realise the open affine subset
$H|_U\!\times_U\! X|_U$
as the subvariety~$\variety{G(t,x)}$
of~$H|_U\!\times X|_U$;
decomposing the ideal~$(G(t,x))$
will therefore give us a model for~$R$.

\begin{lemma}
\label{lemma:square-root-of-s}
	Let~$s$ be the polynomial in~$\FF_{q}[t]$
	defined by
	\begin{equation}
	\label{eq:s-definition}
		s := 
			\begin{array}[t]{l}
			f_0^3 - f_0^2f_1g_2 - 2f_0^2f_2g_1 + f_0^2f_2g_2^2 
			+ f_0f_1^2g_1 + 3f_0f_1f_2g_0 - f_0f_1f_2g_1g_2 
		\\ {}
			- 2f_0f_2^2g_0g_2 + f_0f_2^2g_1^2 - f_1^3g_0 
			+ f_1^2f_2g_0g_2 - f_1f_2^2g_0g_1 + f_2^3g_0^2
		,
		\end{array}
	\end{equation}
	and let~$\alpha$ be its leading coefficient.
	Then~$s$ has a square root in~$\FF_{q}(\sqrt{\alpha})[t]$.
\end{lemma}
\begin{proof}
	The polynomial~$s$ 
	is a square in~$\FF_{q}(\sqrt{\alpha})[t]$
	if and only if
	each of its roots in~$\FFbar_{q}$
	occur with multiplicity~$2$.
	In the notation of~\eqref{eq:preimage-of-Q-on-X},
	we have
	\[
		s(t(Q)) = F(x(P_1))F(x(P_2))F(x(P_3)) ,
	\]
	so~$s(t(Q)) = 0$ if and only if~$F(x(P_i)) = 0$
	for some~$1 \le i \le 3$
	--- that is, if and only if at least one of the~$P_i$ 
	is a Weierstrass point of~$H$.
	But the trigonal map~$g$ was constructed precisely
	so that the Weierstrass points of~$H$
	appear in pairs in the fibres of~$g$:
	hence exactly two of the~$P_i$ must be Weierstrass points,
	and so~$F(x(P_1))F(x(P_2))F(x(P_3)) = 0$
	and~$s(t(Q)) = 0$
	with multiplicity~$2$.
\qed
\end{proof}

\begin{proposition}
\label{proposition:b22-t-equation}
	Let~$s$ be the polynomial of Lemma~\ref{lemma:square-root-of-s},
	and let~$\delta_0$,~$\delta_1$,~$\delta_2$, and~$\delta_4$
	be the polynomials in~$\FF_{q^2}[t]$ defined by
	\[
		\begin{array}{r@{\;:=\;}l}
		\delta_4 
		&
		-27g_0^2 + 18g_0g_1g_2 - 4g_0g_2^3 - 4g_1^3 + g_1^2g_2^2
		,
		\\
		\delta_2 
		&
		12f_0g_1 - 4f_0g_2^2 - 18f_1g_0 + 2f_1g_1g_2 + 12f_2g_0g_2 - 4f_2g_1^2
		,
		\\
		\delta_1
		&
		8\sqrt{s},
		\ \text{ and }
		\\
		\delta_0 
		& 
		-4f_0f_2 + f_1^2
		.
		\end{array}
	\]
	On the curve~$X|_U$, we have
	\begin{equation}
	\label{eq:b22-equation}
		\big( 
			\delta_4(t) b_{22}^2 + \delta_2(t) b_{22} + \delta_0(t) 
		\big)^2 
		- 
		\delta_1(t)^2 b_{22} 
		= 
		0
		.
	\end{equation}
\end{proposition}
\begin{proof}
	Consider again the fibre of~$f: X \to \PP{1}$
	over the generic point~$Q = (t)$ of~$U$
	(as in~\eqref{eq:preimage-of-Q-on-X}).
	If~$\{ P_1 + P_2 + P_3, \iota(P_1)+\iota(P_2)+\iota(P_3) \}$
	is a pair of divisors corresponding to one of the points in the fibre,
	then 
	by
	the Lagrange interpolation formula
	the value of~$b_{22}$ 
	at the corresponding point of~${\widetilde{X}}$~is
	\begin{equation}
	\label{eq:b22-interpolation}
		b_{22}
		=
		\left(\sum y(P_i)/((x(P_i) - x(P_j))(x(P_i) - x(P_k)))\right)^2
		,
	\end{equation}
	where the sum is taken over the cyclic permutations~$(i,j,k)$
	of~$(1,2,3)$.
	After interpolating
	for each pair of divisors
	in the fibre,
	an elementary but involved symbolic calculation shows that
	$b_{22}$ satisfies
	\begin{equation}
	\label{eq:d-relation}
		\Big(
			\Delta b_{22}^2 
			- 2\big(\sum_i\Gamma_i\big) b_{22} 
			+ \frac{1}{\Delta}\Big(
				2\big(\sum_i\Gamma_i^2\big) 
				- \big(\sum_i \Gamma_i\big)^2
			\Big)
		\Big)^2
		-
		64 \big(\prod_i\Gamma_i\big)b_{22}
		= 0 ,
	\end{equation}
	where
	\[
		\Gamma_i 
		:= 
		\big(
			f_2(t)x(P_i)^2 
			+ f_1(t)x(P_i)
			+ f_0(t)
		\big)\Delta_i
		= 
		F(x(P_i))\Delta_i
	\]
	with
	\[
		\Delta_{i} := (x(P_j) - x(P_k))^2
	\]
	for each cyclic permutation~$(i,j,k)$ of~$(1,2,3)$,
	and where
	\(
		\Delta
		:= \Delta_1\Delta_2\Delta_3
	\).

	Now~$\Delta$,~$\sum_i \Gamma_i$,~$\sum_i \Gamma_i^2$, 
	and~$\prod_i \Gamma_i$
	are symmetric functions
	with respect to permutations of 
	the points in the fibre~$g^{-1}(Q) = g^{-1}((t))$.
	They are 
	therefore 
	polynomials in the homogeneous elementary symmetric functions
	\[
		e_1 = \sum_i x(P_i),
		\quad 
		e_2 = \sum_{i<j} x(P_i)x(P_j),
		\quad \mbox{and}\quad 
		e_3 = \prod_i x(P_i)
		,
	\]
	which are polynomials in~$t$.
	Indeed, the~$e_i$ are given by the coefficients of
	$G(t,x)$:
	\[
		e_1 = -g_2(t),\quad 
		e_2 = g_1(t),\quad 
		\mbox{ and }\quad 
		e_3 = -g_0(t) .
	\]
	Expressing 
	$\Delta$,~$\sum_i \Gamma_i$,~$\sum_i \Gamma_i^2$, and~$\prod_i \Gamma_i$
	in terms of~$f_0$,~$f_1$,~$f_2$,~$g_0$,~$g_1$, and~$g_2$,
	and substituting the resulting expressions into~\eqref{eq:d-relation},
	we obtain \eqref{eq:b22-equation}.
\qed
\end{proof}

Equation~\eqref{eq:b22-equation}
gives us a (singular) affine plane model for~$X$.
We can also use~\eqref{eq:b22-equation}
to compute a square root for~$b_{22}$ on~$X|_U$:
we have
\[
	b_{22} = \rho^2,
	\quad \text{where}\quad 
	\rho := 
	\frac{
		\delta_4(t) b_{22}^2 + \delta_2(t) b_{22} + \delta_0(t)
	}{
		\delta_1(t)
	}.
\]
Returning to~\eqref{eq:b22-interpolation},
we observe that
$b_{22}$ is a unit on~$X|_U$,
since its zeroes and poles occur only at points~$Q$
where~$g\circ\pi$ is ramified over~$f(Q)$,
and these points were excluded from~$U$.
Since~$\rho$ is the square root of~$b_{22}$,
it must also be a unit on~$X|_U$.

Given a point~$(t,b_{00},\ldots,b_{22})$ of~$X|_U$,
the corresponding pair of divisors of degree $3$ on~$H$
is cut out by the pair of ideals 
\[
	\Big\{
		\Big( 
			G(t,x), 
			y \pm \big(\frac{b_{02}}{\rho} + \frac{b_{12}}{\rho}x + \frac{b_{22}}{\rho}x^2\big) \Big)
	\Big\}
	.
\]
This is precisely the decomposition of~$(G(t,x))$
that we need to compute the correspondence
from~$H|_U$ to~$X|_U$:
we have
$
        \variety{G(t,x)}
	=
	R \cup R'
$,
where
\begin{equation}
\label{eq:correspondence-eqns}
	R =
	\variety{
		G(t,x), 
		y - \frac{1}{\rho} (b_{02} + b_{12}x + b_{22}x^2)
	}
\end{equation}
and
\[
	R'
	=
	\variety{
		G(t,x), 
		y + \frac{1}{\rho}(b_{02} + b_{12}x + b_{22}x^2)
	}
	.
\]
On the level of divisor classes,
the isogeny~$\phi: \Jac{H} \to \Jac{X}$ is made explicit by the map
\[
	\phi = (\pi_X)_*\circ(\pi_H)^*  ,
\]
where
$\pi_H: R \to H$ and~$\pi_X: R \to X|_U$
are the natural projections
defined by
\(
	(x,y,t,b_{00},\ldots,b_{22}) \mapsto (x,y)
\)
and
\(
	(x,y,t,b_{00},\ldots,b_{22}) \mapsto (t,b_{00},\ldots,b_{22})
\), respectively.
In terms of ideals cutting out effective divisors,
$\phi$ is realized by the map 
$$
	I_D 
	\longmapsto
	\left(I_D + \Big(G(t,x), y - \frac{1}{\rho}\big(b_{02} + b_{12}x + b_{22}x^2\big)\Big)\right)
	\cap \FF_{q}[s,t,b_{00},\ldots,b_{22}]
	.
$$
Taking~$R'$
in place of~$R$ in the above
gives an isogeny equal to~$-\phi$.
It remains to determine the field of definition of~$\phi$.

\begin{proposition}
\label{proposition:X-rationality-criterion}
	If~$S$ is a subgroup in~$\tractables{H}$
	with an~$\FF_{q}$-rational trigonal map~$g$ defined over~$\FF_{q}$,
	and~$s(t)$ is the polynomial defined in Lemma~\ref{lemma:square-root-of-s}.
	then the explicit trigonal construction on~$g$
	described above
	yields an isogeny defined over~$\FF_{q}$
	if and only if
	the leading coefficient of~$s(t)$ is a square in~$\FF_{q}$.
\end{proposition}
\begin{proof}
	We noted earlier
	that~$X|_U$ is defined over the field of definition of~$g$.
	The correspondence~$R$,
	and hence the induced isogeny~$\phi$,
	are both defined over 
	the field of definition of~$\rho$,
	which is the field of definition 
	of $\delta_4\delta_1$, $\delta_2/\delta_1$, and $\delta_0/\delta_1$.
	But $\delta_4$, $\delta_4$, and $\delta_0$
	are all defined over $\FF_{q}$
	(cf.~Proposition~\ref{proposition:b22-t-equation}),
	while $\delta_1$ is defined over $\FF_{q}(\sqrt{\alpha})$
	where $\alpha$ is the leading coefficient of $s$
	by Lemma~\ref{lemma:square-root-of-s}.
\qed
\end{proof}

\begin{remark}
	If~$\phi$ is not defined over~$\FF_{q}$,
	then the Jacobian~$\Jac{X}$ is in fact a quadratic twist 
	of the quotient~$\Jac{H}/S$
	(see~\S\ref{section:other-isogenies}).
	In fact, when~$\phi$ is not~$\FF_{q}$-rational,
	Frobenius exchanges~$\rho$ and~$-\rho$,
	hence~$R$ and~$R'$,
	and therefore
	$\phi$ and~$-\phi$.
	This is a concrete realization
	of the Galois cohomology referred to 
	in the proof of Proposition~\ref{proposition:kernels} below:
	the obstruction to the existence of an isomorphism
	from~$\Jac{H}/S$ to~$\Jac{X}$ over~$\FF_{q}$
	is in fact the interaction of~$\Galois$
	with~$[\pm 1]$ on~$\Jac{X}$.
\end{remark}

If we assume that
the leading coefficients of the polynomials~$s(t)$ 
are uniformly distributed 
for randomly chosen~$H$,~$S$, and~$g$,
then the probability that
$s$ is a square in~$\FF_{q}[t]$
is~$1/2$.
Indeed, it is easily seen that
$s(t)$ is a square for~$H$
if and only if it is not a square for
the quadratic twist of~$H$. 
Suppose~$H: w^2 = \homog{F}(u,v)$
is a hyperelliptic curve.
Let~$c$ be a non-square in~$\FF_{q}$,
and let~$H': w^2 = c\homog{F}(u,v)$
be the quadratic twist of~$H$.
Suppose~$S$ in~$\tractables{H}$
is a tractable subgroup,
represented by a set~$\{ F_1, F_2, F_3, F_4 \}$
of quadratic factors of~$\homog{F}$.
The set~$\{ cF_1, F_2, F_3, F_4 \}$
is a factorization of~$c\homog{F}$,
so it represents a tractable subgroup~$S'$ in~$\tractables{H'}$.
We noted in~\S\ref{section:computing-trigonal-maps}
that scalar multiples of quadratic polynomials
do not affect the construction of trigonal maps;
so if~$S$ has a trigonal map~$g$ defined over~$\FF_{q}$,
then~$g$ is also a trigonal map for~$S'$.
Let~$s$ be the polynomial computed from~$g$ and~$S$
in Lemma~\ref{lemma:square-root-of-s},
and let~$s'$ be the corresponding polynomial
computed for~$g$ and~$S'$.
Looking at the form of~\eqref{eq:s-definition},
we see that~$s'(t) = c^3s(t)$.
Therefore, the leading coefficient of~$s'$
is a square if and only if the leading coefficient of~$s$
is \emph{not} a square.
In particular,
if~$S$ has a trigonal map defined over~$\FF_{q}$,
then so does~$S'$,
and we can construct an isogeny of Jacobians with kernel~$S$
if and only if we cannot construct an isogeny of Jacobians with kernel~$S'$.

This suggests that
the probability that 
we can compute an isogeny defined over~$\FF_{q}$
given a randomly chosen~$H$ and~$S$ in~$\tractables{H}$
with a trigonal map defined over~$\FF_{q}$
is~$1/2$
--- since we have a~$50\%$ chance of being on the ``right'' quadratic twist of~$H$.
This hypothesis is consistent with our experimental observations.
\begin{hypothesis}
\label{hypothesis:rationality}
	For a randomly chosen hyperelliptic curve~$H$ 
	and a uniformly randomly chosen subgroup~$S$ in~$\tractables{H}$
	with a trigonal map~$g$ defined over~$\FF_{q}$,
	the probability that we can compute
	an $\FF_{q}$-rational isogeny~$\phi$ 
	with kernel~$S$ is~$1/2$.
\end{hypothesis}

\section{Computing Isogenies}
\label{section:examples}

Now we will put the ideas above into practice.
Suppose we are given a hyperelliptic curve~$H$ 
of genus~$3$ over~$\FF_{q}$,
and a DLP instance in~$\Jac{H}(\FF_{q})$ to solve.
Our goal is to compute a nonsingular plane quartic curve~$C$ 
and an explcit isogeny~$\Jac{H} \to \Jac{C}$
defined over~$\FF_{q}$,
so that we can solve our DLP instance in~$\Jac{C}(\FF_{q})$.

We begin by computing the set~$\tractables{H}$
of~$\FF_{q}$-rational tractable subgroups of the~$2$-torsion subgroup~$\Jac{H}[2](\FFbar_{q})$
(see Appendix~\ref{appendix:computing-SH} below).
For each~$S$ in~$\tractables{H}$,
we apply Proposition~\ref{proposition:trigonal-map-rationality-criterion}
to determine whether there exists an~$\FF_{q}$-rational
trigonal map~$g$ for~$S$.
If so,
we use the formulae of~\S\ref{section:computing-trigonal-maps}
to compute~$g$;
if not,
we move on to the next~$S$.
Having computed~$g$,
we apply Proposition~\ref{proposition:X-rationality-criterion}
to determine whether we can compute an isogeny over~$\FF_{q}$.
If so,
we use the formulae of~\S\ref{section:equations-for-the-isogeny}
to compute equations for~$X$
and the isogeny~$\phi: \Jac{H} \to \Jac{X}$;
if not,
we move on to the next~$S$.

The formulae of~\S\ref{section:equations-for-the-isogeny}
give an affine model of~$X$ in~$\AA{1}\times\AA{6}$.
In order to apply Diem's algorithm
to the DLP in~$\Jac{X}$,
we need a nonsingular plane quartic model of~$X$:
that is, a nonsingular curve~$C \subset \PP{2}$ isomorphic to~$X$,
cut out by a quartic form.
Such a model exists if and only if~$X$ is not hyperelliptic.
To find~$C$,
we compute a basis~$\mathcal{B} = \{\psi_1,\psi_2,\psi_3\}$ 
of the Riemann--Roch space
of a canonical divisor of~$X$.
This is a routine geometrical calculation;
Hess~\cite{Hess} describes an efficient approach.
In practice,
the algorithms implemented in Magma~\cite{Magma--language,Magma} 
compute~$\mathcal{B}$ very quickly.
The three functions in~$\mathcal{B}$
define a map~$\psi: X \to \PP{2}$,
mapping~$P$ to~$(\psi_1(P):\psi_2(P):\psi_3(P))$.
Up to automorphisms of~$\PP{2}$,
the map~$\psi$ is independent of the choice of basis~$\mathcal{B}$,
and depends only on~$X$.
If the image of~$\psi$ is a conic
(that is, if the~$\psi_i$ satisfy a quadratic relation),
then~$X$ is hyperelliptic;
in this situation we move on to the next~$S$,
since we will gain no advantage from index calculus on~$X$.
Otherwise,
the image of~$\psi$ is a nonsingular plane quartic~$C$,
and~$\psi$ restricts to an isomorphism~$\psi: X \to C$.

If the procedure outlined above succeeds for some~$S$ in~$\tractables{H}$,
then
we have computed
an explicit~$\FF_{q}$-rational isogeny~$\psi_*\circ\phi: \Jac{H} \to \Jac{C}$.
We can then
map our DLP from~$\Jac{H}(\FF_q)$ into~$\Jac{C}(\FF_q)$,
and solve it using Diem's algorithm.

We emphasize that the entire procedure is very fast:
the curve~$X$ and the isogeny can be constructed using
just a few low-degree polynomial operations and some low-dimensional
linear algebra
(and hence
the procedure is polynomial-time in $\log q$, the size of the base field).
For a rough idea of the computational effort involved,
given a random~$H$ over a~$160$-bit prime field
with a tractable subgroup~$S$ in~$\tractables{H}$,
a na\"ive implementation of our algorithms in 
Magma 
computes the trigonal map~$g$,
the curve~$X$, the nonsingular plane quartic~$C$, 
and the isogeny~$\phi: \Jac{H} \to \Jac{C}$
in a few seconds on a 1.2GHz laptop.
Since the difficulty of the construction
depends only upon the difficulty of arithmetic in~$\FF_{q}$
(and \emph{not} upon the size of the DLP subgroup of~$\Jac{H}(\FF_{q})$),
we may conclude that 
instances of the DLP in~$160$-bit Jacobians chosen for cryptography
may also be reduced to instances of the DLP in non-hyperelliptic Jacobians 
in very little time.

\begin{example}
	We will give an example over a small field.
	Let~$H$ be the hyperelliptic curve over~$\FF_{37}$
	defined by
	$$
		H: 
		y^2 
		= 
		x^7 + 28 x^6 + 15 x^5 + 20 x^4 + 33 x^3 + 12 x^2 + 29 x + 2
		.
	$$
	Using the ideas in \S\ref{section:kernel},
	or the algorithms in Appendix~\ref{appendix:computing-SH},
	we find that
	$\Jac{H}$ has one~$\FF_{37}$-rational
	tractable subgroup:
	\[
		\tractables{H}
		=
		\left\{ S \right\}
		\text{\quad where\quad }
		S = 
			\left\{
			\begin{array}{l}
  				u^2 + \xi_1 uv + \xi_2 v^2,
				\ 
  				u^2 + \xi_1^{37} uv + \xi_2^{37} v^2,
			\\
  				u^2 + \xi_1^{37^2} uv + \xi_2^{37^2} v^2 ,
				\ 
				uv + 20v^2
			\end{array}
			\right\}
		,
	\]
	where~$\xi_1$ is an element of~$\FF_{37^3}$
	satisfying~$\xi_1^3 + 29\xi_1^2 + 9\xi_1 + 13 = 0$,
	and~$\xi_2 = \xi_1^{50100}$.
	Applying the methods of~\S\ref{section:computing-trigonal-maps},
	we compute 
	a trigonal map~$g: x \mapsto N(x)/D(x)$ for~$S$,
	taking
	\[
		N(x) = x^3 + 16 x + 22
		\quad \mbox{and}\quad 
  		D(x) = x^2 + 32 x + 18 
		;
	\]
	clearly~$g$
	is 
	defined over~$\FF_{37}$.
	The formulae of~\S\ref{section:equations-for-the-isogeny}
	give us a curve~$X \subset \AA{1}\times\AA{6}$ of genus~$3$,
	defined by
	\[
		X = \variety{
		\begin{array}{l}
			\scriptstyle
			(18t^2 + 15t)b_{22} + (36t + 30)b_{12} + b_{00} + 19t^5 + 10t^4 + 12t^3 + 7t^2 + t + 30,
			\\
			\scriptstyle 
			(32t^2 + 2t + 15)b_{22} + (27t + 5)b_{12} + 2b_{01} + 5t^5 + 26t^4 + 15t^3 + 23t^2 + 19t + 17,
			\\
			\scriptstyle 
			(t^2 + 32t + 21)b_{22} + 2tb_{12} + 2b_{02} + b_{11} + 36t^5 + 29t^4 + 7t^3 + 13t^2 + 21t + 18
			,
			\\
			\scriptstyle 
			b_{00}b_{11} - b_{01}^2,
			b_{00}b_{12} - b_{01}b_{02},
			b_{00}b_{22} - b_{02}^2,
			b_{02}b_{11} - b_{01}b_{12} ,
			b_{02}b_{12} - b_{01}b_{22} ,
			b_{12}^2 - b_{11}b_{22} 
		\end{array}
		}
		.
	\]
	The map on divisors inducing an isogeny from 
	$\Jac{H}$ to~$\Jac{X}$
	with kernel~$S$
	is
	induced by the correspondence~$R$ defined as in~\eqref{eq:correspondence-eqns}
	with
	\[
		\begin{array}{r@{\;}l}
			G(t,x) = &
			x^3 - t x^2 - (32 t - 16) x -18 t + 22
			,
			\\
			\delta_0 = &
			27t^{10} + 20t^9 + 33t^8 + 6t^7 + 16t^6 + 8t^5 + 9t^4 + 2t^3 + 31t^2 + 15t + 16
			,
			\\
			\delta_1 = &
			35t^3 + 8t^2 + 33t + 3
			,
			\\
			\delta_2 = &
20t^7 + 18t^6 + 29t^5 + 14t^4 + 6t^3 + 20t^2 + 12t + 16
			,
			\quad \text{and}
			\\
			\delta_4 = &
				27t^4 + 36t^3 + 13t^2 + 21t
			.
			\\
		\end{array}
	\]
	Computing the canonical morphism of~$X$,
	we find that~$X$ is non-hyperelliptic,
	and isomorphic to the nonsingular plane quartic
	curve
	\[
		C
		= 
		\variety{
			\begin{array}{c}
			u^4 + 26u^3v + 2u^3w + 17u^2v^2 + 9u^2vw + 20u^2w^2 + 34uv^3 + 24uv^2w 
			\\ {} 
			+ 5uvw^2 + 36uw^3 + 19v^4 + 13v^3w + v^2w^2 + 23vw^3 + 5w^4

			\end{array}
		}
		.
	\]
	Composing the isomorphism with the isogeny~$\Jac{H} \to \Jac{X}$,
	we obtain an explicit isogeny
	$\phi: \Jac{H} \to \Jac{C}$.
	We can verify that~$\Jac{H}$ and~$\Jac{C}$ are isogenous
	by checking that the zeta functions of~$H$ and~$C$ are identical:
	indeed, direct calculation with Magma shows that
	\[
		Z(H;T) 
		= 
		Z(C;T) 
		= 
		\frac{
			37^3 T^6 + 4\cdot37^2 T^5 - 6\cdot37 T^4 - 240 T^3 
			- 6 T^2 + 4 T + 1
		}{
			(37T - 1)(T - 1)
		}
		.
	\]
	Let~$D =  [ (10:28:1) - (14:6:1) ]$ 
	and~$D' = [ (19:28:1) - (36:13:1) ]$ 
	be divisor classes on~$H$;
	we have
	$D' = [22359]D$.
	Applying~$\phi$,
	we find that
	\[
	\begin{array}{r@{\;=\;}l}
		\phi(D) 
		&
  		\left[ 
			(7 : 18 : 1) + (34 : 34 : 1)
			- (18 : 22 : 1) - (15 : 33 : 1)
		\right]
		\text{\quad  and}\\
		\phi(D') 
		&
  		\left[ 
  			(7 : 23 : 1) + (6 : 13 : 1)
  			- (13 : 15 : 1) - (7 : 18 : 1)
		\right]
		;
	\end{array}
	\]
	direct calculation verifies that~$\phi(D') = [22359]\phi(D)$,
	as expected.
\end{example}

\section{Expectation of Existence of Computable Isogenies}
\label{section:probabilities}

Our aim in this section
is
to estimate the proportion of genus~$3$ hyperelliptic Jacobians
over~$\FF_{q}$
for which the methods of this article
produce an~$\FF_{q}$-rational isogeny
--- and thus for which the DLP may be solved using Diem's algorithm ---
as~$q$ tends to infinity.
We will assume that 
if we are given a selection of~$\FF_q$-rational tractable subgroups
of a given Jacobian,
then the probabilities that each 
will yield a rational isogeny are mutually independent.
This hypothesis appears to be consistent with our experimental observations.
\begin{hypothesis}
\label{hypothesis:independence}
	For a randomly chosen hyperelliptic curve $H$,
	the probabilities that 
	we can compute an $\FF_{q}$-rational isogeny
	with kernel $S$ for each $S$ in $\tractables{H}$
	are mutually independent.
\end{hypothesis}

%
%
%

\begin{theorem}
\label{theorem:success-probability}
	Assume Hypotheses~\ref{hypothesis:hyperellipticity},
	\ref{hypothesis:trigonal-map-rationality},
	\ref{hypothesis:rationality},
	and~\ref{hypothesis:independence}.
	As~$q$ tends to infinity,
	the expectation that 
	the algorithms in this article
	will give a reduction of the DLP 
	in a subgroup of~$\Jac{H}(\FF_{q})$
	for a randomly chosen hyperelliptic curve~$H$
	of genus~$3$ over~$\FF_{q}$
	to a subgroup of~$\Jac{C}(\FF_{q})$
	for some nonsingular plane quartic curve~$C$
	is
	\begin{equation}	
	\label{eq:expectation}
		\sum_{T \in \mathcal{T}} 
		\Big(
			\big(1 - (1 - 1/4)^{s(T)}\big)
			/
			\prod_{n \in T}
			\big(\nu_T(n)! \cdot n^{\nu_T(n)}\big) 
		\Big)
		\approx
		0.1857 ,
	\end{equation}
	where
	$\mathcal{T}$
	denotes the set of integer partitions of~$8$
	and 
	$\nu_T(n)$
	denotes the multiplicity of an integer~$n$ in a partition~$T$,
	and~$s(T) = \#\tractables{H}$,
	where~$H$ is \emph{any} hyperelliptic curve over~$\FF_{q}$
	such that the multiset of degrees of the~$\FF_{q}$-irreducible factors
	of its hyperelliptic polynomial coincides with~$T$.
\end{theorem}
\begin{proof}
	Suppose~$H$ is a randomly chosen hyperelliptic curve of genus~$3$
	over~$\FF_{q}$.
	Hypotheses~\ref{hypothesis:hyperellipticity}, 
	\ref{hypothesis:trigonal-map-rationality},
	and \ref{hypothesis:rationality}
	together
	imply that 
	for each~$S$ in~$\tractables{H}$,
	the probability that we can compute an isogeny with kernel~$S$
	defined over~$\FF_{q}$
	is
	\( 1/2\cdot1/2\cdot1 = 1/4 \).
	Hypothesis~\ref{hypothesis:independence}
	implies that we have an equal chance of constructing an isogeny
	from each~$S$ in~$\tractables{H}$,
	so the probability that we can compute an isogeny
	over~$\FF_{q}$
	from~$\Jac{H}$
	is
	\( 1 - (1 - 1/4)^{\#\tractables{H}} \).
	The expectation that we can compute an isogeny over~$\FF_{q}$
	given a curve over~$\FF_{q}$
	is therefore 
	\begin{equation} 
	\label{eq:success-probability}
		E_q := 
		\frac{
			\sum_{\homog{F}} ( 1 - (3/4)^{\#\tractables{H}} ) 
		}{
			\sum_{\homog{F}}  1
		},
	\end{equation}
	where~$H$ is the curve defined by~$w^2 = \homog{F}(u,v)$,
	and~$\homog{F}$ ranges over
	the set of all homogeneous squarefree polynomials 
	of degree~$8$ over~$\FF_{q}$.
	Lemma~\ref{lemma:number-of-tractable-subgroups}
	implies that~$\#\tractables{H}$
	depends only on the degrees of 
	the~$\FF_{q}$-irreducible factors of~$\homog{F}$,
	so the map~$T \mapsto s(T)$ is well-defined. 
	For each~$T$ in~$\mathcal{T}$,
	let~$ N_q(T)~$
	denote the number of homogeneous squarefree polynomials over~$\FF_{q}$
	whose multiset of degrees of~$\FF_{q}$-irreducible factors
	coincides with~$T$.
	We can now rewrite~\eqref{eq:success-probability}
	as
	\[
		E_q
		= 
		\frac{
			\sum_{T \in \mathcal{T}}(1 - (3/4)^{s(T)})N_q(T)
		}{
			\sum_{T\in \mathcal{T}}N_q(T)
		}
		.
	\]
	There are
	$N_q(n) = \frac{1}{n}\sum_{d|n}\mu(d)q^{n/d}$
	monic irreducible 
	polynomials of degree~$n$
	over~$\FF_{q}$
	(here $\mu$ is the M\"obius function).
	Clearly
	\(
		N_q(T) 
		= 
		(q - 1)\prod_{n \in T}\binom{N_q(n)}{\nu_T(n)} 
	\),
	so
	\[
		N_q(T) 
		= 
		\Big(\prod_{n\in T}(\nu_T(n)!\cdot n^{\nu_T(n)})\Big)^{-1}q^9 + O(q^8) 
		,
	\]
	and 
	$\sum_{T \in \mathcal{T}}N_q(T) = q^9 + O(q^8)$.
	Therefore,
	as~$q$ tends to infinity,
	we have
	$$
		\lim_{q \to \infty} E_q
		=
		\sum_{T \in \mathcal{T}} 
		\Big(
			\left(1 - (3/4)^{s(T)}\right)
			/
			\prod_{n \in T}
			\left(\nu_T(n)! \cdot n^{\nu_T(n)}\right) 
		\Big)
		.
	$$
	The result follows upon
	explicitly computing this sum,
	using the values for~$s(T)$ listed in 
	Lemma~\ref{lemma:number-of-tractable-subgroups}.
\qed
\end{proof}

Theorem~\ref{theorem:success-probability}
gives the expectation of our ability to construct
an explicit isogeny for
a randomly selected hyperelliptic curve.
However,
looking at the table in Lemma~\ref{lemma:number-of-tractable-subgroups},
we see that we can be sure that a particular
curve has no isogenies with tractable kernels
defined over~$\FF_{q}$
if we use only curves whose hyperelliptic polynomials have
an irreducible factor of degree~$5$ or~$7$
(or a single irreducible factor of degree~$3$).
It may be difficult to efficiently construct a curve in this form
if we are using a CM construction,
for example,
to ensure that the Jacobian has a large prime-order subgroup.
In any case,
it is interesting to note that
the security of genus~$3$ hyperelliptic Jacobians
depends significantly upon the factorization of their hyperelliptic polynomials.
This observation has no analogue for elliptic curves
or Jacobians of curves of genus~$2$.
Of course,
if~$E: y^2 = F(x)$ is an elliptic curve
and~$F$ is completely reducible,
then~$\#E(\FF_{q})$ is divisible by~$4$,
and in particular~$\#E(\FF_{q})$ cannot be prime;
but this does not reduce the security of~$E(\FF_{q})$
to the extent that a completely reducible hyperelliptic polynomial
does for a curve of genus~$3$.

\begin{remark}
\label{remark:trigonal-map-descent}
	We noted in \S\ref{section:trigonal-construction}
	that the~$\FFbar_{q}$-isomorphism class of the curve~$X$
	in the trigonal construction is independent
	of the choice of trigonal map.
	If there is no trigonal map defined over~$\FF_{q}$ 
	for a given subgroup~$S$ in~$\tractables{H}$,
	then the methods of~\S\ref{section:computing-trigonal-maps}
	construct a pair of Galois-conjugate trigonal maps~$g_1$ and~$g_2$
	(corresponding to the roots of~\eqref{eq:lambda-eqn})
	instead.
	Applying the trigonal construction to~$g_1$ and~$g_2$,
	we obtain curves~$X_1$ and~$X_2$ over~$\FF_{q^2}$.
	If the isomorphism between~$X_1$ and~$X_2$
	were made explicit,
	then we could descend it to compute
	a curve~$X$ over~$\FF_{q}$ 
	in the~$\FFbar_{q}$-isomorphism class of~$X_1$ and~$X_2$,
	and hence a nonsingular plane quartic~$C$ 
	over~$\FF_{q}$
	and an isogeny~$\Jac{H} \to \Jac{C}$.
	We note that the isogeny may not be defined over~$\FF_{q}$,
	but this approach 
	could still allow us to replace the~$1/4$ in~\eqref{eq:success-probability}
	and~\eqref{eq:expectation}
	with~$1/2$,
	raising the expectation of success
	in Theorem~\ref{theorem:success-probability}
	to~$31.13\%$.
\end{remark}

\begin{example}
	Let $p = 1008945029102471339$.  Note that $p$ is a 60-bit prime;
	if $H$ is a hyperelliptic curve of genus $3$ over $\FF_{p}$
	such that $\Jac{H}(\FF_{p})$ has a large prime-order subgroup
	and 
	if Gaudry--Thom\'e--Th\'eriault--Diem index calculus is the
	fastest algorithm for solving DLP instances in $\Jac{H}(\FF_{p})$,
	then 
	$\Jac{H}$ has roughly the same security level as an elliptic curve
	over a 160-bit field.

	We generated one million random hyperelliptic curves of genus~$3$
	over~$\FF_{p}$ using Magma.
	For each curve $H$,
	we computed the set $\tractables{H}$ 
	of tractable subgroups;
	then,
	for each $S$ in $\tractables{H}$
	we determined whether there was 
	an $\FF_{p}$-rational trigonal map for $S$,
	and if so
	whether there was an $\FF_{p}$-rational isogeny
	with kernel~$S$.
	Of these curves, $502005$ (that is, $50.02\%$)
	had at least one rational tractable subgroup.
	Between them, the $10^6$ curves had $1002244$ rational tractable
	subgroups,
	of which
	$501629$ had a rational trigonal map
	(that is, $50.05\%$, which is close to the $50\%$ predicted by 
	Hypothesis~\ref{hypothesis:trigonal-map-rationality}).
	Of these subgroups,
	$250560$ led to a rational isogeny
	(that is, $49.95\%$,
	which is close to the $50\%$ predicted by 
	Hypothesis~\ref{hypothesis:rationality}).
	We found that $185814$ of the curves
	had at least one $\FF_{p}$-rational isogeny,
	none of which had a hyperelliptic codomain
	(this is compatible with 
	Hypothesis~\ref{hypothesis:hyperellipticity}).
	In particular,
	we could move a discrete logarithm problem
	for $18.58\%$ of these curves
	(recall that Theorem~\ref{theorem:success-probability}
	predicts a success rate of about $18.57\%$).

\end{example}

\section{Other Isogenies}
\label{section:other-isogenies}

So far,
we have concentrated on using isogenies
with kernels generated by differences of Weierstrass points
to move instances of the DLP
from hyperelliptic to non-hyperelliptic Jacobians.
More generally,
we could use isogenies with other kernels.
There are two important issues to consider here:
the first is 
a theoretical restriction on the types of
subgroups that can be kernels
of isogenies of Jacobians,
and the second is a practical restriction
on the isogenies that we can currently compute.

Let~$H$ be a hyperelliptic curve of genus~$3$.
We want to characterize the subgroups~$S$ of~$\Jac{H}$
that are kernels of isogenies of Jacobians,
combining standard results
from the theory of abelian varieties with some special
results on curves of genus~$3$.
For our purposes,
it is enough to know that the~$l$-\emph{Weil pairing} is a 
nondegenerate, bilinear pairing on the~$l$-torsion 
of an abelian variety,
which can be efficiently evaluated in the case where the abelian variety
is the Jacobian of a hyperelliptic curve;
for further detail, we refer the reader to~\cite[Ex.~A.7.8]{Hindry--Silverman}.

\begin{definition}
	Let~$A$ be an abelian variety over $\FF_{q}$,
	and let~$l$ be a positive integer
	coprime with $q$.
	We say a subgroup~$S$ of~$A[l]$
	is \emph{maximal~$l$-isotropic}
	if
	\begin{enumerate}
		\item	the~$l$-Weil pairing
			on~$A[l]$ restricts trivially to~$S$, and
		\item	$S$ is not properly contained in any other 
			subgroup of~$A[l]$ satisfying~(1).
	\end{enumerate}
\end{definition}

If~$l$ is a prime not dividing $q$,
then every maximal~$l$-isotropic subgroup of~$\Jac{H}(\FFbar_{q})[l]$
is isomorphic to~$(\ZZ/l\ZZ)^3$.
The situation is more complicated when~$l$ is not prime:
for example,~$\Jac{H}[2]$ is a maximal~$4$-isotropic subgroup of~$\Jac{H}[4]$,
but it is isomorphic to $(\ZZ/2\ZZ)^6$ and not $(\ZZ/4\ZZ)^3$.

\begin{proposition}
\label{proposition:kernels}
	Let~$H$ be a hyperelliptic curve of genus~$3$
	over~$\FF_{q}$
	such that~$\Jac{H}$ is absolutely simple.
	Let~$S$ be a finite, nontrivial,~$\FF_{q}$-rational 
	subgroup of~$\Jac{H}(\FFbar_{q})$.
	There exists a curve~$X$ of genus~$3$
	over~$\FF_{q}$,
	and an isogeny~$\phi: \Jac{H} \to \Jac{X}$
	with kernel~$S$,
	if and only if
	$S$ is a maximal~$l$-isotropic subgroup of~$\Jac{H}[l]$
	for some positive integer~$l$.
	The isogeny~$\phi$ is defined over~$\FF_{q^2}$.
\end{proposition}
\begin{proof}
	The quotient~$\Jac{H} \to \Jac{H}/S$
	always exists as an isogeny of abelian varieties,
	and is defined over~$\FF_{q}$
	(see Serre~\cite[\S III.3.12]{Serre-AC}).
	For the quotient 
	to be an isogeny of Jacobians,
	there must be an integer~$l$
	such that~$S$ is a maximal~$l$-isotropic subgroup
	(see Proposition~16.8 of Milne~\cite{Milne}):
	this ensures that 
	the canonical polarization on~$\Jac{H}$
	induces a principal polarization 
	on the quotient~$\Jac{H}/S$.
	The theorem of Oort and Ueno~\cite{Oort--Ueno}
	therefore guarantees that there will be an isomorphism
	of principally polarized abelian varieties
	over~$\FFbar_{q}$
	from~$\Jac{H}/S$ to the Jacobian~$\Jac{X}$ of some
	irreducible curve $X$
	(irreducibility of~$X$
	follows from the fact that~$\Jac{H}$, and hence~$\Jac{H}/S$,
	is absolutely simple).
	Composing this isomorphism
	with the quotient map
	gives an isogeny of Jacobians from $\Jac{H}$ to $\Jac{X}$
	with kernel $S$.
	Standard arguments from Galois cohomology
	(see Serre~\cite[\S III.1]{Serre-GC}, for example)
	show that the isomorphism is defined over either~$\FF_{q}$ or~$\FF_{q^2}$,
	and it follows that the isogeny $\Jac{H} \to \Jac{X}$
	must be defined over $\FF_{q}$ or $\FF_{q^2}$.
\qed
\end{proof}

\begin{remark}
	Proposition~\ref{proposition:kernels}
	does \emph{not} hold in higher genus:
	for every~$g \ge 4$,
	there are~$g$-dimensional abelian varieties
	that are not isomorphic to Jacobians.
	Indeed, this is the generic situation:
	for $g \ge 2$
	the moduli space of $g$-dimensional abelian varieties
	is~$g(g+1)/2$-dimensional,
	with the Jacobians occupying 
	a subspace of dimension~$(3g - 3)$
	--- which is strictly less than $g(g+1)/2$ for $g \ge 4$.
	We should not therefore expect an arbitrary quotient of a Jacobian
	to be isomorphic to a Jacobian in genus $g \ge 4$.
	Proposition~\ref{proposition:kernels}
	does hold in genus~$1$ and~$2$,
	and in these cases
	the isogenies are always defined over~$\FF_{q}$.
\end{remark}

We can expect the curve~$X$ of Proposition~\ref{proposition:kernels}
to be non-hyperelliptic.
To compute an~$\FF_{q}$-rational isogeny from~$\Jac{H}$
to a non-hyperelliptic Jacobian, therefore,
the minimum requirement is 
an~$\FF_{q}$-rational~$l$-isotropic subgroup 
of~$\Jac{H}(\FFbar_{q})$
isomorphic to~$(\ZZ/l\ZZ)^3$ for some prime~$l$.
We emphasize
that this subgroup 
need \emph{not} be contained in~$\Jac{H}(\FF_{q})$.
Indeed,
there may be isogenies from~$\Jac{H}$ to non-hyperelliptic Jacobians
over~$\FF_{q}$
even when~$\Jac{H}(\FF_{q})$ has prime order
(which would be the desirable situation in cryptological applications).

The major obstruction to using more general isogenies
to move DLP instances 
is the lack of general constructions for explicit isogenies
in genus~$3$.
Apart from integer multiplications, automorphisms, Frobenius isogenies,
and the construction for isogenies with tractable kernels
exhibited above,
we know of no constructions for explicit isogenies
of general Jacobians of genus~$3$ hyperelliptic curves.
In particular,
while we know that the curve $X$ of Proposition~\ref{proposition:kernels}
exists, we generally have no means of computing a defining equation for it,
let alone equations for a correspondence between $H$ and $X$
that would allow us to move DLP instances from $\Jac{H}$ to $\Jac{X}$.
This situation stands in marked contrast
to the case of isogenies of elliptic curves,
which have been made completely explicit by V\'elu~\cite{Velu}.
Deriving general formulae 
for explicit isogenies in genus~$3$ (and~$2$)
remains a significant problem in computational number theory.

\subsection*{Acknowledgements}
\label{section:acknowledgements}

The greater part of this work was completed in the Department of Mathematics 
at Royal Holloway, University of London,
where the author was supported by EPSRC grant EP/C014839/1.
The author gratefully acknowledges Roger Oyono and Christophe Ritzenthaler
for discussions which inspired this research,
and Steven Galbraith 
and the anonymous referees
for their helpful suggestions.


\appendix

\section{Appendix: Computing~$\tractables{H}$}
\label{appendix:computing-SH}

Given a hyperelliptic curve $H$ of genus $3$ over $\FF_{q}$,
we want to compute the set~$\tractables{H}$
of~$\FF_{q}$-rational tractable subgroups of~$\Jac{H}$.
Algorithm~\ref{algorithm:subgroup-enumeration}
splits the hyperelliptic polynomial of $H$ into Galois orbits of factors,
before calling the recursive subroutine 
Algorithm~\ref{algorithm:subgroup-subroutine}
to enumerate $\tractables{H}$.
This algorithm is 
included only for completeness,
and is 
not particularly efficient
(we suggest some optimisations 
in Remark~\ref{remark:tractable-implementation} 
below.)

\begin{algorithm}
\label{algorithm:subgroup-enumeration}
	Given a hyperelliptic curve~$H$ of genus~$3$ 
	over~$\FF_{q}$,
	enumerates the set~$\tractables{H}$
	of~$\FF_{q}$-rational tractable subgroups of~$\Jac{H}[2](\FFbar_{q})$.
	Each subgroup in~$\tractables{H}$
	is represented by
	a set of four coprime quadratic factors of~$\homog{F}$.
	\begin{description}
		\item[Input] 
			The hyperelliptic polynomial~$\homog{F}(u,v)$ of~$H$.
		\item[Output] 
			The set~$\tractables{H}$.
		\item[Step 1]
			Let~$\mathcal{F}$
			be the set of irreducible factors of~$\homog{F}$ 
			over its splitting field,
			\\
			scaled so that~$\homog{F} = \prod_{L\in\mathcal{F}}L$,
			and set~$\mathcal{O} := \{ \}$.
		\item[Step 2]
			Choose a polynomial~$L$ from~$\mathcal{F}$.
			Set~$O := (L)$,
			set~$\mathcal{F} := \mathcal{F}\setminus\{L\}$,
			\\
			and 
			set~$L_1 := L$.
		\item[Step 3]
			Set~$L := \sigma(L)$,
			where~$\sigma$ denotes the~$q^\mathrm{th}$ power Frobenius map.
			\\
			If~$L \not= L_1$,
			then append~$L$ to~$O$,
			set~$\mathcal{F} := \mathcal{F}\setminus\{L\}$,
			and go to \textbf{Step~3}.
			\\
			If~$L = L_1$,
			then set~$\mathcal{O} := \mathcal{O}\cup\{ O \}$;
			if~$\mathcal{F} \not= \emptyset$,
			then go to \textbf{Step~2}.
		\item[Step 4]
			Return the result of 
			Algorithm~\ref{algorithm:subgroup-subroutine}
			applied to~$\mathcal{O}$.
	\end{description}
\end{algorithm}

\begin{algorithm}
\label{algorithm:subgroup-subroutine}
	Given a set of~$\Galois$-orbits of coprime linear polynomials 
	over $\FFbar_{q}$,
	returns the~$\Galois$-invariant sets of coprime quadratic products of the polynomials.
\begin{description}
	\item[Input]
		A set~$\mathcal{O}$
		of disjoint sequences
		of distinct linear polynomials.
		Each sequence~$O = (O_1,\ldots,O_m)$ 
		in~$\mathcal{O}$
		must satisfy~$O_1 = \sigma(O_{m})$ 
		and~$O_{i+1} = \sigma(O_i)$
		for~$1 \le i < m$,
		where~$\sigma$ denotes the~$q^\mathrm{th}$-power
		Frobenius map.
	\item[Output]
		The set~$\mathcal{S}$ 
		of~$\Galois$-stable sets of coprime quadratic polynomials
		such that
		$
			\prod_{S \in \mathcal{S}}\prod_{Q \in S} Q
			=
			\prod_{O \in \mathcal{O}}\prod_{L \in O} L
		$.
	\item[Step 1]
		If~$\mathcal{O}$ is empty,
		then return~$\mathcal{S} := \{ \emptyset \}$.
	\item[Step 2]
		Choose a sequence~$O$ from~$\mathcal{O}$,
		and set~$m := \#O$.
		\\
		If~$m$ is even,
		then 
		let~$\mathcal{T}$ be 
		the result of Algorithm~\ref{algorithm:subgroup-subroutine}
		applied to~$\mathcal{O}\setminus\{O\}$,
		\\
		and set
		$
			\mathcal{S} 
			:= 
			\{ \{ O_{i}\cdot O_{(m/2)+i}: 1\le i \le m/2 \} \cup T : T \in \mathcal{T} \} 
		$.
		\\
		If~$m$ is odd, 
		then set~$\mathcal{S} := \{ \}$.
	\item[Step 3]
		For each~$P$ in~$\mathcal{O} \setminus \{ O \}$
		such that~$\#P = \#O = m$,
		\begin{description}
			\item[Step 3i]
				Set 
				$
					\mathcal{U} 
					:= 
					\left\{\!\{ O_{1+i}\cdot P_{1+((i+j)\bmod m)} : 0\le i < m \} : 0 \le j < m \right\}
				$.
			\item[Step 3ii]
				Let~$\mathcal{V}$
				be the result of 
				Algorithm~\ref{algorithm:subgroup-subroutine}
				applied to~$\mathcal{O}\setminus\{O,P\}$.
			\item[Step 3iii]
				Set
				$
					\mathcal{S}
					:=
					\mathcal{S}
					\cup
					\left\{
						U \cup V : 
						U \in \mathcal{U},
						\ 
						V \in \mathcal{V}
					\right\}
				$.
		\end{description}
	\item[Step 4]
			Return~$\mathcal{S}$.
\end{description}
\end{algorithm}

\begin{remark}
\label{remark:tractable-implementation}
	As we noted above,
	Algorithms $4$ and $5$ are not particularly efficient:
	for conceptual simplicity we worked over the splitting field
	of the hyperelliptic polynomial,
	and this can be extremely slow in practice.
	A number of simple optimizations 
	will significantly improve the performance of this algorithm:
	the key is to avoid field extensions where possible,
	and to minimize their degree in any case.
	Before factoring~$\homog{F}$ over its splitting field
	we should factor it over~$\FF_{q}$,
	and then work on a case-by-case basis
	depending on the degrees of the~$\FF_{q}$-irreducible factors.
	For example,
	if~$\homog{F}$ has an odd number of odd-degree factors,
	then $\tractables{H}$ is empty
	by Lemma~\ref{lemma:number-of-tractable-subgroups},
	and we can simply return the empty set.
	If~$\homog{F}$ is~$\FF_{q}$-irreducible,
	then it is not necessary to factor $\homog{F}$ over its splitting field
	(which is~$\FF_{q^8}$):
	there is one tractable subgroup,
	and it corresponds to the four quadratic factors of $\homog{F}$
	that we obtain by factoring $\homog{F}$ over $\FF_{q^4}$.
	Making similar modifications for the cases where
	$\homog{F}$ has factors of degree $6$,
	we can avoid working over any extensions 
	of degree greater than $4$.
	If desired, we can further avoid some field extensions
	in the case where $\homog{F}$ has only low-degree factors.
	These modifications resulted in a factor-of-$50$ speedup
	for our experiments with $60$-bit prime fields;
	the unmodified Algorithms~\ref{algorithm:subgroup-enumeration}
	and~\ref{algorithm:subgroup-subroutine}
	should \emph{not} be used in practice.
\end{remark}

\end{document}